\documentclass[12pt,amscd]{amsart}
\usepackage{amssymb}
\numberwithin{equation}{section}
\newtheorem{thm}[equation]{Theorem}

\newtheorem{lem}[equation]{Lemma}

\newenvironment{pf}{\proof[\proofname]}{\endproof}
\newenvironment{pf*}[1]{\proof[#1]}{\endproof}

\theoremstyle{definition}

\newtheorem{defn}[equation]{Definition}

\theoremstyle{remark}

\newtheorem*{rmk}{Remark}

\newtheorem*{ack}{Acknowledgement}

\newcommand{\comment}[1]{}

\makeatother

\begin{document}
\baselineskip=18truept


\def\C {{\mathbb C}}
\def\Cn {{\mathbb C}^n}
\def\R {{\mathbb R}}
\def\Rn {{\mathbb R}^n}
\def\Z {{\mathbb Z}}
\def\N {{\mathbb Z}_{>0}}
\def\cal#1{{\mathcal #1}}
\def\bb#1{{\mathbb #1}}

\def\dbar {\bar \partial }
\def\dir {{\mathcal D}}
\def\lev#1{{\mathcal L}\left(#1\right)}
\def\lap {\Delta }
\def\ol {{\mathcal O}}
\def\E {{\mathcal E}}
\def\J {{\mathcal J}}
\def\U {{\mathcal U}}
\def\V {{\mathcal V}}
\def\z {\zeta }
\def\Harm {\text {Harm}\, }
\def\grad {\nabla }
\def\dexh {\{ M_k \} _{k=0}^{\infty } }
\def\sing#1{#1_{\text {sing}}}
\def\reg#1{#1_{\text {reg}}}
\def\pione#1{\pi_1(#1)}
\def\pioneat#1#2{\pi_1(#1,#2)}

\def\lquotient{\big{\backslash}}

\def\setof#1#2{\{ \, #1 \mid #2 \, \} }

\def\image#1{\text{\rm im}\,\bigl[#1\bigr]}
\def\kernel#1{\text{\rm ker}\,\bigl[#1\bigr]}

\def\holecl {M\setminus \overline M_0}
\def\hole {M\setminus M_0}

\def\nd{\frac {\partial }{\partial\nu } }
\def\ndof#1{\frac {\partial#1}{\partial\nu } }

\def\pdof#1#2{\frac {\partial#1}{\partial#2}}
\def\pdwrt#1{\frac{\partial}{\partial #1}}

\def\cinf{\ensuremath{\mathcal C^{\infty}} }
\def\cinfns{\ensuremath{\mathcal C^{\infty}}}

\def\dist{\text{\rm dist}}

\def\diam{\text{\rm diam}}

\def\real{\text{\rm Re}\, }

\def\imag{\text{\rm Im}\, }

\def\supp{\text{\rm supp}\, }

\def\Vol{\text{\rm vol}}

\def\restrict#1{{\upharpoonright_{#1}}}

\def\sm{\setminus}

\def\plshclass{{\mathcal {P}}}
\def\strplshclass{{\mathcal S\mathcal P}}

\def\geqtrace#1#2{{\geq_{(#1,#2)}}}
\def\gtrace#1#2{{>_{(#1,#2)}}}
\def\leqtrace#1#2{{\leq_{(#1,#2)}}}
\def\ltrace#1#2{{<_{(#1,#2)}}}



\def\anal{analytic }
\def\analns{analytic}

\def\bdd{bounded }
\def\bddns{bounded}

\def\cpt{compact }
\def\cptns{compact}

\def\cpx{complex }
\def\cpxns{complex}

\def\cont{continuous }
\def\contns{continuous}

\def\dime{dimension }
\def\dimens{dimension }

\def\exh{exhaustion }
\def\exhns{exhaustion}

\def\fn{function }
\def\fnns{function}

\def\fns{functions }
\def\fnsns{functions}

\def\holo{holomorphic }
\def\holons{holomorphic}

\def\mero{meromorphic }
\def\merons{meromorphic}

\def\holoconvex{holomorphically convex }
\def\holoconvexns{holomorphically convex}

\def\ircomp{irreducible component }
\def\concomp{connected component }
\def\ircompns{irreducible component}
\def\concompns{connected component}
\def\ircomps{irreducible components }
\def\concomps{connected components }
\def\ircompsns{irreducible components}
\def\concompsns{connected components}

\def\irred{irreducible }
\def\irredns{irreducible}

\def\con{connected }
\def\conns{connected}

\def\comp{component }
\def\compns{component}
\def\comps{components }
\def\compsns{components}

\def\mfld{manifold }
\def\mfldns{manifold}
\def\mflds{manifolds }
\def\mfldsns{manifolds}

\def\nbd{neighborhood }
\def\nbds{neighborhoods }
\def\nbdns{neighborhood}
\def\nbdsns{neighborhoods}

\def\harm{harmonic }
\def\harmns{harmonic}
\def\plh{pluriharmonic }
\def\plhns{pluriharmonic}
\def\plsh{plurisubharmonic }
\def\plshns{plurisubharmonic}

\def\qplsh#1{$#1$-plurisubharmonic}
\def\hplsh{$(n-1)$-plurisubharmonic }
\def\hplshns{$(n-1)$-plurisubharmonic}

\def\para{parabolic }
\def\parans{parabolic}

\def\rel{relatively }
\def\relns{relatively}

\def\str{strictly }
\def\strns{strictly}

\def\strg{strongly }
\def\strgns{strongly}

\def\cvx{convex }
\def\cvxns{convex}

\def\wrt{with respect to }
\def\wrtns{with respect to}

\def\st {such that }
\def\stns {such that}

\def\hm {harmonic measure }
\def\hmns {harmonic measure}

\def\hmib {harmonic measure of the ideal boundary of }
\def\hmibns {harmonic measure of the ideal boundary of}

\def\Vert{\text{{\rm Vert}}\, }
\def\Edge{\text{{\rm Edge}}\, }

\def\til#1{\tilde{#1}}
\def\wtil#1{\widetilde{#1}}

\def\what#1{\widehat{#1}}

\def\seq#1#2{\{#1_{#2}\} }


\def\vphi {\varphi }


\def\inv{   ^{-1}  }

\def\ssp#1{^{(#1)}}

\def\set#1{\{ #1 \}}

\title[A cup product lemma]
{A cup product lemma for continuous plurisubharmonic
functions}
\author[T.~Napier]{Terrence Napier}
\address{Department of Mathematics\\Lehigh University\\Bethlehem, PA 18015\\USA}
\email{tjn2@lehigh.edu}
\thanks{Preprint of an article submitted for consideration in \emph{Journal of Topology and Analysis},  \copyright\,  2016 World Scientific Publishing Company, www.worldscientific.com/worldscinet/jta}
\author[M.~Ramachandran]{Mohan Ramachandran}
\address{Department of Mathematics\\University at Buffalo\\Buffalo, NY 14260\\USA}
\email{ramac-m@buffalo.edu}


\date{October 16, 2016}

\begin{abstract}
A version of Gromov's cup product lemma 
in which one factor is the $(1,0)$-part of the differential of
a \cont \plsh \fn is obtained. As an application, it
is shown that a \con non\cpt complete K\"ahler manifold that 
has exactly one end and admits a \cont \plsh \fn that
is \str \plsh along some 
germ of a $2$-dimensional \cpx \anal set at some point
has the Bochner--Hartogs property; that is,
the first compactly supported cohomology with values in the
structure sheaf vanishes.
\end{abstract}

\maketitle

\section*{Introduction} \label{introduction}

Versions of Gromov's cup product lemma
have been applied in various settings in order
to obtain results concerning the \holo structure of
a \con non\cpt complete K\"ahler manifold $(X,g)$. 
For example, in \cite{Gro-Sur la groupe fond}, \cite{Li Structure
complete Kahler}, \cite{Gro-Kahler hyperbolicity},
\cite{Gromov-Schoen}, \cite{NR-Structure theorems},
\cite{Delzant-Gromov Cuts}, \cite{NR Filtered ends}, and \cite{NR L2 Castelnuovo}, versions of the
cup product lemma, along with analogues
of the classical Castelnuovo--de Franchis theorem,
yield conditions under which $X$ admits a proper
\holo mapping onto a Riemann surface. 
We may illustrate the general approach to such results as follows.  According to Theorem~0.2 of~\cite{NR L2 Castelnuovo}, which may be viewed as a version of the
cup product lemma, $\theta_1\wedge\theta_2\equiv 0$ for any pair of closed \holo $1$-forms $\theta_1$ and $\theta_2$
on $X$ \st $\theta_1$ is bounded, $\real\theta_1$ is exact, and $\theta_2$ is in $L^2$. If these $1$-forms are
linearly independent and $X$ has bounded geometry, then 
Stein factorization of the mapping $\theta_1/\theta_2$ of~$X$ 
into~$\mathbb P^1$ gives a proper \holo mapping onto a Riemann surface (see
Theorem~0.1 and Corollary~0.3 of \cite{NR L2 Castelnuovo}).  

Corollary~5.4 of~\cite{NR L2 Castelnuovo} is a version of the above cup product lemma in which the factor $\theta_1$ is 
replaced with $\partial\vphi$, where $\vphi$
is a \cinf \plsh \fn with bounded gradient. One can remove the boundedness
condition on the gradient of~$\vphi$ by instead requiring~$\theta_2$ to be in~$L^2$ with respect to the
(complete) K\"ahler metric 
\(h\equiv g+\lev{e^{2\vphi}}\)  (we denote the Levi form of a \fnns~$\psi$ by $\lev\psi$), since $e^\vphi$ has bounded gradient
with respect to~$h$. The main goal of the present paper is to obtain
analogues in which the \plsh \fnns~$\vphi$ is allowed to be only \contns. Not surprisingly,
for the proofs, one applies the above construction of~$h$ to terms of 
a sequence of \cinf approximations of~$\vphi$ (in place of~$\vphi$) provided by
a theorem of Greene and Wu (see Corollary 2 of
Theorem~4.1 of~\cite{Greene-Wu} as well as \cite{Richberg} and
\cite{Dem}).
The following is the simplest form of these analogues:
\begin{thm}\label{cup product cont very simple version intro thm}
 Let $\vphi$ be a \cont \plsh \fn on a \con non\cpt complete K\"ahler
manifold $(X,g)$.  Then there exists a (complete)
K\"ahler metric~$h$ on~$X$ \st $h\geq g$ and \st
 \(\partial\vphi\wedge\theta\equiv 0\) as a current for every (closed) \holo $1$-form~$\theta$
on~$X$ that is in~$L^2$ with respect to~$h$.
\end{thm}

Theorem~\ref{cup product cont very simple version intro thm} is applied in this paper to obtain a result
which is related to those of \cite{Ramachandran-BH for coverings},
 \cite{NR-BH Weakly 1-complete}, \cite{NR-BH Regular hyperbolic Kahler}, and \cite{NR BH Bounded Geometry},
 and which we now describe.
A \con non\cpt \cpx manifold~$X$ for which $H^1_c(X,\ol)=0$ is said to
have the {\it Bochner--Hartogs property} (see
Hartogs~\cite{Hartogs}, Bochner~\cite{Bochner}, and
Harvey~and~Lawson~\cite{Harvey-Lawson}). Equivalently, for every
\cinf compactly supported form~$\alpha$ of type $(0,1)$ with
$\dbar\alpha=0$ on~$X$, there is a \cinf compactly supported \fn
$\beta$ on~$X$ \st $\dbar\beta=\alpha$.  A \con non\cpt \cpx manifold~$X$ with the Bochner--Hartogs
property must satisfy certain analytic and topological conditions.  For example, $X$ must have exactly one end,
every \holo \fn on the \con complement of a \cpt subset of~$X$ must extend holomorphically to~$X$, and
$X$ cannot admit a proper \holo mapping onto a Riemann surface. 
Some further elementary consequences are described in~\cite{NR BH Bounded Geometry}. Examples of manifolds
of dimension~$n>1$ having the Bochner--Hartogs property include
strongly $(n-1)$-complete \cpx manifolds (Andreotti and
Vesentini~\cite{Andreotti-Vesentini}) and strongly
hyper-$(n-1)$-convex K\"ahler manifolds (Grauert and
Riemenschneider~\cite{Grauert-Riemenschneider}). 
According to Theorem~3.3 of \cite{NR  BH Bounded Geometry}, a one-ended
\con non\cpt complete K\"ahler manifold~$X$ that is hyperbolic (in the sense that $X$ admits a positive symmetric Green's \fnsns)
and has no nontrivial $L^2$ \holo $1$-forms has the Bochner--Hartogs property. 
This observation and Theorem~\ref{cup product cont very simple version intro thm} together lead to the following (cf.~Proposition~4.4 of~\cite{NR Weak Lefschetz}):
\begin{thm}\label{BH plsh fn cont two pos ev thm}
Let $(X,g)$ be a \con non\cpt complete K\"ahler manifold with 
exactly one end. Assume
that $X$ admits a \cont \plsh \fnns~$\vphi$
whose restriction to some $2$-dimensional germ of an \anal set 
at some point $p\in X$ is \str \plsh  
(for example, a $\cal C^2$ \plsh \fn on~$X$ whose Levi form
has at least two positive eigenvalues at some point has this property).  Then
$H^1_c(X,\ol)=0$.
\end{thm}
\begin{pf*}{Sketch of the proof}
We may assume that $\vphi\geq 0$. Given a $\dbar$-closed compactly supported \cinf
$(0,1)$-form~$\alpha$ on~$X$, we may fix a nonempty \con \cpt set~$K$, a positive
constant $a<\sup\vphi$, and a \concompns~$\Omega$ of
$\setof{x\in X}{\vphi(x)<a}$ \st $\supp\alpha\subset K\subset\Omega$, $\Omega\sm K$ is \conns, and $p\in K$. 
Nakano~\cite{Nk}, Greene and Wu~\cite{Greene-Wu}, and
Demailly~\cite{Dem} provide a complete K\"ahler metric~$h$ on~$\Omega$ with respect
to which~$\Omega$ is hyperbolic, and we may choose~$h$ as in
Theorem~\ref{cup product cont very simple version intro thm}. The condition on~$\vphi$ then implies
that $(\Omega,h)$ has no nontrivial
$L^2$ \holo $1$-forms, and
Theorem~3.3 of \cite{NR  BH Bounded Geometry} then gives the claim.
\end{pf*}
Further applications of
the versions of the cup product lemma appearing in
this paper will be considered elsewhere (see~\cite{NR Weakly Special}).

In Section~\ref{ends basic sect}, we
recall some terminology and basic facts concerning ends. In Section~\ref{greens fn sect}, we recall some terminology and facts from
potential theory, as well as
some notation concerning Hermitian metrics. 
In Section~\ref{cinf approx plsh sect}, we
recall some facts concerning \cinf approximation
of \plsh \fnsns. 
Section~\ref{cup product cont sect} contains
the statements and proofs of the desired versions of the cup product lemma,
of which Theorem~\ref{cup product cont very simple version intro thm} is a special case.
Section~\ref{bh plsh two pos ev sect} contains the proof of
Theorem~\ref{BH plsh fn cont two pos ev thm}, which is
obtained as special case of a version that allows for
multiple ends.

\begin{ack}  
The authors would like to thank Cezar Joita for some useful
references.
\end{ack}

\section{Ends}\label{ends basic sect}

In this section we
recall some terminology and basic facts concerning ends.
By an {\it end} of a \con manifold~$M$, we will mean either a \comp $E$
of $M\setminus K$ with non\cpt closure, where $K$ is a given \cpt
subset of $M$, or an element of
\[
\lim _{\leftarrow} \pi _0 (M\setminus K),
\]
where the limit is taken as $K$ ranges over the \cpt subsets of
$M$ (or equivalently, the \cpt subsets of $M$ for which the complement $M\setminus K$
has no \rel \cpt \compsns, since the union of any \cpt subset of~$M$ with the \rel \cpt \concomps
of its complement is \cptns). The number of ends of $M$ will be
denoted by~$e(M)$. For a \cpt set $K$ \st $M\setminus K$ has no
\rel \cpt \compsns, we will call
\[
M\setminus K=E_1\cup \cdots \cup E_m,
\]
where $\seq Ej_{j=1}^m$ are the distinct \comps of $M\setminus K$,
an {\it ends decomposition} for~$M$. The following elementary lemma
will allow us to modify ends decompositions and to pass to 
ends decompositions in domains:
\begin{lem}\label{basic ends lem}
Let $M$ be a \con non\cpt \cinf manifold.
\begin{enumerate}

\item[(a)] Given an ends decomposition $M\sm K=E_1\cup\cdots\cup
E_m$, there is a \con \cpt set~$K'\supset K$ \st any
domain~$\Theta$ in~$M$ containing~$K'$ has an ends decomposition
$\Theta\sm K=E_1'\cup\cdots\cup E_m'$, where $E_j'=E_j\cap\Theta$
for $j=1,\dots,m$.

\item[(b)] If $\Omega$ and $\Theta$ are domains in~$M$ with
$\Theta\subset\Omega$, and both $M\sm\Omega$ and $\Omega\sm\Theta$
have no \cpt \compsns, then $M\sm\Theta$ has no \cpt \compsns.

\item[(c)] If $E$ is an end of~$M$, then there exists an end $A_0$ \st
$\overline A_0\subset E$ and $E\sm A_0\Subset M$.

\item[(d)] If $E$ is an end of~$M$, $F_1,\dots,F_k\subset E$
are disjoint ends of~$M$ for which $k>1$ and $E\sm(F_1\cup\cdots\cup F_k)\Subset M$, and 
$s\in\set{2,\dots,k}$, then there exist disjoint ends 
$A_1,\dots,A_s\subset E$ of~$M$ \st
$\overline A_j\subset F_j$ and $F_j\sm A_j\Subset M$ for $j=1,\dots,s-1$, and $A_s\supset F_s\cup\cdots\cup F_k$.

\end{enumerate}
\end{lem}
\begin{pf}
Parts (a) and (b) are provided by Lemma~1.1 of \cite{NR BH Bounded Geometry}. We  will prove
part~(d), and leave to the reader 
the proof of part~(c) (which is easier and mostly contained in the proof of part~(d)). 
The claim is trivial in
dimension~$1$, so we may assume that $n\equiv\dim M>1$. We may fix a nonempty 
\cinf domain~$\Omega_0$ \st $\overline E\sm(F_1\cup\cdots\cup F_k)\subset\Omega_0\Subset M$, $M\sm\Omega_0$ has no
\cpt \concompsns, and $E\cap\Omega_0$, as well as $F_j\cap\Omega_0$ for $j=1,\dots,k$, are \conns. Suppose $G_1$ and $G_2$ are two
distinct \comps of $F_1\sm\overline\Omega_0$ and therefore, of $E\sm\overline\Omega_0$ and of $M\sm\overline\Omega_0$.
Fixing a boundary \comp $B_i$ of $G_i$ for $i=1,2$, we see that $B_1$ and $B_2$ are distinct
boundary \comps of~$\Omega_0$, and there exists a \con \cpt set $C\subset (F_1\cap\Omega_0)\cup B_1\cup B_2$
\st $C\cap B_1$ and $C\cap B_2$ are singletons, and the sets $\Omega_0\sm C$, $E\cap\Omega_0\sm C$, and
$F_j\cap\Omega_0\sm C$ for $j=1,\dots,k$, are \conns, and are therefore ends of $\Omega_0\sm C$.
For example, we may take $C$ to be the image of 
an embedded \cinf path in~$M$ \st the intial point lies in $B_1$, the terminal point lies in~$B_2$, all other
points lie in $F_1\cap\Omega_0$, and the path meets $B_1$ and $B_2$ transversely. 
In dimension $n>1$, such a path is locally nonseparating, while in dimension $n=1$, $B_1$ and $B_2$ admit
disjoint \nbds in~$F_1$  \st the intersection of each of these \nbds with $\Omega_0$ is a coordinate annulus
whose intersection with $C$ is a radial line segment. 
One \concomp of $M\sm(\Omega_0\sm C)=(M\sm\Omega_0)\cup C$ meets, and therefore
contains, 
the \con non\cpt set $\overline G_1\cup\overline G_2\cup C$.  The remaining \concomps 
are the (finitely many)
\concomps of $M\sm\Omega_0$ which do not meet this set, and are in fact, precisely the 
\concomps of
$(M\sm\Omega_0)\sm(\overline G_1\cup\overline G_2)$. 
Thus each \concomp of $M\sm(\Omega_0\sm C)$ is non\cpt and $F_1\sm(\Omega_0\sm C)$
has \str fewer \comps than $F_1\sm\Omega_0$ (and $M\sm(\Omega_0\sm C)$ has \str fewer
\comps than $M\sm\Omega_0$).  We may fix a \cinf domain~$\Omega_1$ \st
\[
\overline E\sm(F_1\cup\cdots\cup F_k)\subset\Omega_1\Subset\Omega_0\sm C,
\]
$(\Omega_0\sm C)\sm\Omega_1$
has no \cpt \concompsns, and $(E\cap\Omega_0\sm C)\cap\Omega_1=E\cap\Omega_1$, as well as
$(F_j\cap\Omega_0\sm C)\cap\Omega_1=F_j\cap\Omega_1$ for $j=1,\dots,k$, are ends of $\Omega_1$.
Thus we get a \cinf domain
\(\Omega_2\equiv\Omega_1\cup(\Omega_0\sm F_1)\) \st 
$\overline E\sm(F_1\cup\cdots\cup F_k)\subset\Omega_2\subset\Omega_0\sm C$; $(\Omega_0\sm C)\sm\Omega_2$ and 
$M\sm(\Omega_0\sm C)$, and therefore
$M\sm\Omega_2$, have no \cpt \compsns; and $E\cap\Omega_2$, as well as $F_j\cap\Omega_2$ for 
$j=1,\dots,k$, are \conns. Moreover, $F_j\sm\Omega_2=F_j\sm\Omega_0$ for $j=2,\dots,k$,
and each \comp of $F_1\sm\Omega_2=F_1\sm\Omega_1$ contains a \comp of $F_1\sm(\Omega_0\sm C)$,
so $F_1\sm\Omega_2$ has strictly fewer \comps than $F_1\sm\Omega_0$.
Proceedingly inductively, we get a \cinf 
domain~$\Omega$ \st $\overline E\sm(F_1\cup\cdots\cup F_k)\subset\Omega\Subset M$; $E\cap\Omega$ 
and $F_j\cap\Omega$ for $j=1,\dots,k$ are \conns; 
and for each $j=1,\dots,s-1$,
$A_j\equiv F_j\sm\overline\Omega$ is \con and satisfies $\overline A_j\subset F_j$. The ends
$A_1,\dots,A_{s-1},A_s\equiv(E\cap\Omega)\cup F_s\cup\cdots\cup F_k$ then
have the properties listed in part~(d).
\end{pf}

\section{Green's \fns and hermitian metrics}\label{greens fn sect}

In this section we recall some terminology and facts from
potential theory, as well as
some notation concerning Hermitian metrics. 
A \con non\cpt oriented Riemannian manifold~$(M,g)$ is called
\textit{hyperbolic} if there exists a positive symmetric Green's
\fnns~$G(x,y)$ on~$M$; otherwise, $M$ is called \emph{parabolic}.
Equivalently, $M$ is hyperbolic if given a \rel \cpt \cinf
domain~$\Omega $ for which no \concomp of $M\sm\Omega$ is \cptns,
there are a \concompns~$E$ of $M\sm\overline\Omega$ and a (unique)
greatest \cinf \fn $u_E:\overline E\to[0,1)$ \st $u_E$ is \harm on
$E$, $u_E=0$ on $\partial E$, and $\sup_Eu_E=1$.
We will also call $E$, and any end containing~$E$, a \emph{hyperbolic} end.  An end
that is not hyperbolic is called \emph{parabolic}, and 
we set $u_E\equiv 0$ for any parabolic end \compns~$E$ of $M\sm\overline\Omega$. 
We call the \fn $u:M\sm\Omega\to[0,1)$
defined by $u\restrict{\overline E}=u_E$ for each \concompns~$E$
of $M\sm\overline\Omega$, the \emph{\harm measure of the ideal
boundary of}~$M$ \emph{with respect to} $M\sm\overline\Omega$. A
sequence $\seq x\nu$ in $M$ with $x_\nu\to\infty$ and
$G(\cdot,x_\nu)\to 0$ (equivalently, $u(x_\nu)\to 1$) is called a
\emph{regular sequence}. Such a sequence always exists (for $M$
hyperbolic). A sequence $\seq x\nu$ tending to infinity with
$\liminf_{\nu\to\infty}G(\cdot,x_\nu)>0$ (i.e.,
$\limsup_{\nu\to\infty}u(x_\nu)<1$ or equivalently, $\seq x\nu$
has no regular subsequences) is called an \emph{irregular
sequence}. Clearly, every sequence tending to infinity that is not
regular admits an irregular subsequence. We say that an end $E$ of
$M$ is
\emph{regular} (\emph{irregular}) if every sequence in~$E$ tending to
infinity in~$M$ is regular (respectively, there exists an irregular
sequence in~$E$). Another characterization of hyperbolicity is that $M$
is hyperbolic if and only if $M$ admits a nonconstant negative
\cont sub\harm \fnns~$\vphi$. In fact, if $\seq x\nu$ is a
sequence in $M$ with $x_\nu\to\infty$ and $\vphi(x_\nu)\to 0$,
then $\seq x\nu$ is a regular sequence.

The \emph{energy} (or \emph{Dirichlet integral}) of
a suitable \fnns~$\vphi$ (for example, a \fn with first-order
distributional derivatives) on a Riemannian manifold~$M$ is given
by $\int_M|\nabla\vphi|^2\,dV$. As is well known, 
the \hmib an oriented \con non\cpt Riemannian manifold 
has finite energy.

Let $X$ be a \cpx manifold with almost
\cpx structure $J\colon TX\to TX$. By a \emph{Hermitian metric}
on~$X$, we will mean a Riemannian metric~$g$ on $X$ \st
$g(Ju,Jv)=g(u,v)$ for every choice of real tangent vectors $u,v\in
T_pX$ with $p\in X$. We call $(X,g)$ a \emph{Hermitian manifold}.
We will also denote by $g$ the \cpx bilinear extension of $g$ to
the complexified tangent space $(TX)_{\C}$. The corresponding real
$(1,1)$-form~$\omega$ is given by $(u,v)\mapsto\omega(u,v)\equiv
g(Ju,v)$. The corresponding Hermitian metric (in the sense of a
smoothly varying family of Hermitian inner products) in the \holo
tangent bundle $T^{1,0}X$ is given by \((u,v)\mapsto g(u,\bar
v)\). Observe that with this convention, under the \holo vector
bundle isomorphism $(TX,J)\overset\cong\to T^{1,0}X$ given by
$u\mapsto \frac 12(u-iJu)$, the pullback of this Hermitian metric
to~$(TX,J)$ is given by \((u,v)\mapsto\tfrac 12g(u,v)-\tfrac
i2\omega(u,v)\). In a slight abuse of notation, we will also
denote the induced Hermitian metric in $T^{1,0}X$, as well as the
induced Hermitian metric in
$\Lambda^r(TX)_{\C}\otimes\Lambda^s(T^*X)_{\C}$, by~$g$. The
corresponding Laplacians are given by:
\begin{align*}
\lap&=\lap_d\equiv -(dd^*+d^*d),\\
\lap_{\dbar}&=-(\dbar\dbar^*+\dbar^*\dbar),\\
\lap_{\partial}&=-(\partial\partial^*+\partial^*\partial).
\end{align*}
If $(X,g,\omega)$ is \emph{K\"ahler}, i.e.,~$d\omega=0$, then
$\lap=2\lap_{\dbar}=2\lap_{\partial}$.

\section{Smooth approximation of plurisubharmonic functions}\label{cinf approx plsh sect}

For the proof of (the generalizations of) 
Theorem~\ref{cup product cont very simple version intro thm}, we will form suitable \cinf approximations of \cont \plsh \fns 
by applying the following theorem of Greene and Wu (see Corollary 2 of
Theorem~4.1 of~\cite{Greene-Wu} as well as \cite{Richberg}, 
\cite{Dem}, and Lemma~2.15 of~\cite{Demailly-cinf approximation}):
\begin{thm}[Greene--Wu]\label{cinf approximation thm}
Let $(X,g)$ be a Hermitian manifold, let $K\subset X$ be a closed
subset, let $\vphi$ be a \cont \plsh \fn on~$X$ that is of class
\cinf on some \nbd of~$K$, and let $\delta$ be a positive \cont
\fn on~$X$. Then there exists a \cinf \fnns~$\psi$ \st
$\psi=\vphi$ on a \nbd of~$K$, and $\vphi\leq\psi<\vphi+\delta$ and
$\lev\psi\geq-\delta g$ on~$X$.
\end{thm}
That the \fnns~$\psi$ may be chosen to be equal to~$\vphi$ near~$K$
is not included explicitly in the statement in~\cite{Greene-Wu}.  
As in 
the proof of
Lemma~2.15 of~\cite{Demailly-cinf approximation}, one may 
obtain the above version by applying a \cinf version of the maximum appearing in~\cite{Demailly-Cohomology of q-convex spaces}
(where Demailly applied it to obtain a quasi-\plsh \fn with
logarithmic singularities along a given \anal subset). 
The authors have not found the theorem stated in precisely this form in the literature,
so a proof is
provided in this section for the convenience of the reader.
Natural modifications of the proof give
analogous statements for other classes of \fns (see for example, \cite{Joita-Napier-Ramachandran}).  For example, a \cont \str \plsh \fnns~$\vphi$ on a \cpx manifold (or \cpx space)~$X$
that is \cinf on a \nbd of a closed set~$K$ may be approximated on~$X$ by a \cinf \str \plsh \fn that is equal to~$\vphi$ on a \nbd of~$K$.
   
The \cinf maximum
is given by the following lemma, the proof of which is left to the reader:
\begin{lem}[Demailly \cite{Demailly-Cohomology of q-convex spaces}] \label{cinf maximum function lemma}
Let $\kappa\colon\R\to[0,\infty)$ be a \cinf \fn \st
$\supp\kappa\subset (-1,1)$, $\int_{\R}\kappa (u)\,du=1$, and
$\int_{\R}u\kappa (u)\,du=0$. For each $m\in\N$ and each
$r=(r_1,\dots,r_m)\in(\R_{>0})^m$, let $\cal M_r\colon\R^m\to\R$
be the \fn given by
\[
\cal M_r(t)\equiv\int_{\R^m}\biggl[\max_{1\leq j\leq
m}(t_j+r_ju_j)\biggr]\prod_{1\leq j\leq m}\kappa (u_j)\,du_j
\qquad\text{for }t=(t_1,\dots,t_m)\in\R^m.
\]
Then for each $r=(r_1,\dots,r_m)\in(\R_{>0})^m$, $\cal M_r$ has
the following properties:
\begin{enumerate}

\item[(a)] For each $t=(t_1,\dots,t_m)\in\R^m$,
\begin{align*}
\cal M_r(t)&=\int_{\R^m}\biggl[\max_{1\leq j\leq
m}(t_j+u_j)\biggr]\prod_{1\leq j\leq m}r_j\inv\kappa
(u_j/r_j)\,du_j \\
&=\int_{\R^m}\biggl[\max_{1\leq j\leq m}u_j\biggr]\prod_{1\leq
j\leq m}r_j\inv\kappa ((u_j-t_j)/r_j)\,du_j.
\end{align*}

\item[(b)] $\cal M_r$ is \cinf and \cvxns, and for each
$t=(t_1,\dots,t_m)\in\R^m$ and each permutation $\sigma$ of
$\set{1,\dots,m}$, $\cal M_r(t)=\cal
M_{(r_{\sigma(1)},\dots,r_{\sigma(m)})}(t_{\sigma(1)},\dots,t_{\sigma(m)})$.

\item[(c)] For each $j=1,\dots,m$, $0\leq\pdof{\cal M_r}{t_j}\leq 1$.

\item[(d)] For each $s\in\R$, we have $\cal
M_r(t_1+s,\dots,t_m+s)=\cal M_r(t_1,\dots,t_m)+s$.

\item[(e)] For every $t=(t_1,\dots,t_m)\in\R^m$, \(\max_{1\leq
j\leq m}t_j\leq\cal M_r(t)\leq\max_{1\leq j\leq m}(t_j+r_j)\).

\item[(f)] If $t'=(t_0,t_1,\dots,t_m)=(t_0,t)\in\R^{m+1}$ and
$r'=(r_0,r_1,\dots,r_m)\in(\R_{>0})^{m+1}$ with $t_0+r_0\leq
t_1-r_1$, then $\cal M_{r'}(t')=\cal M_r(t)$.

\item[(g)] For $m=1$, $\cal M_r(t)=t$ for each $t\in\R$.

\item[(h)] If $\vphi=(\vphi_1,\dots,\vphi_m)$ is an $m$-tuple of
\cinf real-valued \fns on a \cpx manifold~$X$, then
\[
\lev{\cal M_r(\vphi)}(v,v)\geq\min_{1\leq j\leq m}\lev{\vphi_j}(v,v)\qquad
\forall\,v\in T^{1,0}X.
\]
\end{enumerate}
\end{lem}

\begin{pf*}{Proof of Theorem~\ref{cinf approximation thm}} 
Let $n\equiv\dim X$.  Fixing a \fnns~$\kappa$ as in 
Lemma~\ref{cinf maximum function lemma}, we may form
the corresponding family of \cinf \fns $\set{\cal M_r}$.
We may choose open sets $\seq\Omega
k_{k=0}^3$, $\seq U\nu _{\nu\in N}$, $\seq V\nu_{\nu\in N}$, and
$\seq W\nu_{\nu\in N}$, and \fns $\lambda$ and 
$\seq\alpha\nu_{\nu\in N}$ \st
\begin{enumerate}
\item[(i)] We have
$K\subset\Omega_0\subset\overline\Omega_0\subset\Omega_1
\subset\overline\Omega_1
\subset\Omega_2\subset\overline\Omega_2\subset\Omega_3$;

\item[(ii)] We have $\lambda\in\cinf (X)$, $0\leq\lambda\leq 1$ on
$X$, $\lambda\equiv 1$ on $\Omega_1$, and
$\supp\lambda\subset\Omega_2$;

\item[(iii)] The \fn $\vphi$ is of class \cinf on $\Omega_3$;

\item[(iv)] For each $\nu\in N$, we have $U_\nu\Subset
V_\nu\Subset W_\nu\Subset X\sm K$ and, for $k=1,2,3$,
$W_\nu\Subset\Omega_k$ if
$W_\nu\cap\overline\Omega_{k-1}\neq\emptyset$ and $W_\nu\Subset
X\sm\overline\Omega_{k-1}$ if $W_\nu\cap
X\sm\Omega_k\neq\emptyset$;

\item[(v)] $\seq W\nu$ is locally finite in $X$ and
$\seq U\nu$ covers $X\sm\Omega_0$;

\item[(vi)] For each $\nu\in N$
and each $\epsilon>0$, there exists a \cinf \plsh \fnns~$\rho$
on~$W_\nu$ \st $\vphi\leq\rho\leq\vphi+\epsilon$ on
$W_\nu$  (if $W_\nu\subset\Omega_3$, then we may take
$\rho=\vphi\restrict {W_\nu}$); and

\item[(vii)] For each $\nu\in N$, $\alpha_\nu\in\cinf(W_\nu)$,
$0\leq\alpha_\nu\leq 1$ on $W_\nu$, $\alpha_\nu\equiv 0$ on
$U_\nu$, and $\alpha_\nu\equiv 1$ on $W_\nu\sm V_\nu$.

\end{enumerate}
Given positive constants $\seq\epsilon\nu_{\nu\in N}$  with
$\epsilon_\nu<\delta$ on $W_\nu$ and
\[
\epsilon_\nu|\lev{\alpha_\nu}(v,v)|\leq\delta(x)|v|^2_g\qquad
\forall\,x\in W_\nu,\,v\in T^{1,0}_xX
\]
for each~$\nu$, we may, for each $\nu\in N$, set
\[
r_\nu\equiv\tfrac 14\min \setof{\epsilon_\mu}{\mu\in N,\, W_\mu\cap
W_\nu\neq\emptyset},
\]
and choose a \cinf \plsh \fnns~$\rho_\nu$ on~$W_\nu$ with $\vphi\leq\rho_\nu\leq\vphi+r_\nu$ on~$W_\nu$. If $W_\nu\subset\Omega_3$, then we set
$\rho_\nu\equiv\vphi\restrict {W_\nu}$. Thus we may define a \fn
$\beta$ on $X\sm\overline\Omega_0$ as follows. Given a point $p\in
X\sm\overline\Omega_0$, we let $\nu_1,\dots,\nu_m\in N$ be the
distinct indices with $p\in W_\nu$ if and only if
$\nu\in\set{\nu_1,\dots,\nu_m}$, we set
$r=(r_{\nu_1},\dots,r_{\nu_m})$ and
$\rho=(\rho_{\nu_1}-\epsilon_{\nu_1}\alpha_{\nu_1},\dots,\rho_{\nu_m}-\epsilon_{\nu_m}\alpha_{\nu_m})$,
and we set $\beta(p)\equiv\cal M_r(\rho(p))$. 
We will show that $\beta$ is of class
\cinfns,
$\vphi\leq\beta<\vphi+\delta$, and $\lev\beta\geq-\delta g$
on $X\sm\overline\Omega_0$. We will then show that if the constants $\seq\epsilon\nu$
are sufficiently small, then the \fn
$\psi\equiv\lambda\cdot\vphi+(1-\lambda)\cdot\beta\in\cinf(X)$
satisfies $\lev\psi\geq-\delta g$ at each point in
$\overline\Omega_2\sm\Omega_1$, and it will follow that~$\psi$
has the required properties.

For $p\in X\sm\overline\Omega_0$ and 
$\beta(p)=\cal M_r(\rho(p))$, where
$\nu_1,\dots,\nu_m\in N$,
$r=(r_{\nu_1},\dots,r_{\nu_m})$, and
$\rho=(\rho_{\nu_1}-\epsilon_{\nu_1}\alpha_{\nu_1},\dots,\rho_{\nu_m}-\epsilon_{\nu_m}\alpha_{\nu_m})$ are
as in the previous paragraph,
we may assume that $p\in U_{\nu_1}$, and we may
choose a \rel \cpt \nbdns~$Q$ of~$p$ \st $Q\Subset U_\nu$ ($V_\nu$, $W_\nu$,
$\Omega_k$, $X\sm\overline{U_\nu}$, $X\sm\overline {V_\nu}$,
$X\sm\overline {W_\nu}$, $X\sm\overline \Omega_k$) whenever $p\in
U_\nu$ (respectively, $V_\nu$, $W_\nu$, $\Omega_k$, $X\sm\overline
{U_\nu}$, $X\sm\overline {V_\nu}$, $X\sm\overline {W_\nu}$,
$X\sm\overline \Omega_k$). In particular, by part~(e) of
Lemma~\ref{cinf maximum function lemma}, we have
\begin{align*}
\vphi(p)\leq&\rho_{\nu_1}(p)
=\rho_{\nu_1}(p)-\epsilon_{\nu_1}\alpha_{\nu_1}(p)
\leq
\max_{1\leq j\leq
m}(\rho_{\nu_j}(p)-\epsilon_{\nu_j}\alpha_{\nu_j}(p))
\leq\beta (p)\\
&\leq \max_{1\leq j\leq
m}(\rho_{\nu_j}(p)-\epsilon_{\nu_j}\alpha_{\nu_j}(p)+r_{\nu_j})<\vphi
(p)+\delta (p).
\end{align*}

After reordering, we may assume that for some $k\in\set{1,\dots,m}$, $\nu=\nu_1,\dots,\nu_k\in N$ are precisely those indices
for which $p\in \overline {V_\nu}$ (i.e., for which
$\overline Q\cap\overline {V_\nu}\neq\emptyset$). Setting $s\equiv (r_{\nu_1},\dots,r_{\nu_k})$, we then have
\(\beta=\cal
M_s(\rho_{\nu_1}-\epsilon_{\nu_1}\alpha_{\nu_1},\dots,\rho_{\nu_k}-\epsilon_{\nu_k}\alpha_{\nu_k})\) on
\[
Q\Subset U_{\nu_1}\cap W_{\nu_2}\cap\cdots\cap W_{\nu_m}.
\]
For if $\nu\in
N\sm\set{\nu_1,\dots,\nu_k}$, then on 
\(Q\cap W_\nu\subset
W_\nu\sm\overline {V_\nu}\),
we have
\[
\rho_\nu-\epsilon_\nu\alpha_\nu+r_\nu+r_{\nu_1}\leq
\vphi+r_\nu-\epsilon_\nu+r_\nu+r_{\nu_1}<\vphi\leq\rho_{\nu_1}=\rho_{\nu_1}-\epsilon_{\nu_1}\alpha_{\nu_1}.
\]
Part~(f) of Lemma~\ref{cinf maximum function lemma} now gives the
expression for~$\beta$. Hence~$\beta$ is smooth, and part~(h)
of the lemma implies that $\lev\beta\geq-\delta g$. 

If in the above $p\in\overline\Omega_2\sm\Omega_1$, then 
$W_{\nu_j}\Subset\Omega_3\sm\overline\Omega_0$ and
$\rho_{\nu_j}=\vphi\restrict{W_{\nu_j}}$ for $j=1,\dots,m$, and
hence by part~(d) of  Lemma~\ref{cinf maximum function lemma},
$\beta=\vphi+\xi$ and $\psi=\vphi+(1-\lambda)\xi$ 
on~$Q$, where
$\xi=\cal
M_s(-\epsilon_{\nu_1}\alpha_{\nu_1},\dots,-\epsilon_{\nu_k}\alpha_{\nu_k})$. 
Moreover, for each point
$x\in Q$ and each tangent vector
$v\in T^{1,0}_xX$, we have
\begin{align*}
\lev{(1-\lambda)\xi}(v,v)&=
(1-\lambda(x))\lev\xi(v,v)-2\real
[\partial\lambda(v)\overline{\partial\xi(v)}]
-\lev\lambda(v,v)\xi (x)
\\
&\geq-(1-\lambda(x))\max_{1\leq
j\leq k}\epsilon_{\nu_j}\lev{\alpha_{\nu_j}}(v,v)
\\
&-2
|\partial\lambda(v)|\cdot\biggl(\sum_{j=1}^k\epsilon_{\nu_j}|\partial\alpha_{\nu_j}(v)|\biggr)
-|\lev{\lambda}(v,v)|\max_{1\leq j\leq
k}\epsilon_{\nu_j}.
\end{align*}
Forming a locally finite covering of $\overline\Omega_2\sm\Omega_1$
by such \nbdsns~$Q$, we see that for sufficiently small
$\seq\epsilon\nu$, we have
\[
\lev\psi\geq\lev{(1-\lambda)\xi}\geq-\delta g
\]
at each point in $\overline\Omega_2\sm\Omega_1$, and it follows
that~$\psi$ has the required properties on~$X$.
\end{pf*}

\section{Cup product lemmas}\label{cup product cont sect}

In this section, we consider the promised versions of the cup product lemma, of which
Theorem~\ref{cup product cont very simple version intro thm} is a special case. Theorem~\ref{cup product cont very simple version intro thm} suffices for the applications considered in this paper, but the more general versions are required for
applications to be considered elsewhere (see \cite{NR Weakly Special}).
According to
Theorem~5.2 of~\cite{NR L2 Castelnuovo} (see also Corollary~5.4
of~\cite{NR L2 Castelnuovo}), if $\vphi$ is nonconstant \cinf \plsh \fn
with bounded gradient on a \con complete K\"ahler
manifold~$(X,g)$, then $\partial\vphi\wedge\theta\equiv 0$
for every $L^2$ \holo $1$-form~$\theta$. One may obtain the bounded gradient condition 
by replacing $\vphi$ with $e^\vphi$ and $g$ with $g+\lev{e^{2\vphi}}$. We will
obtain the following version for $\vphi$ \cont (in fact, for a countable family
of \cont \plsh \fnsns):
\begin{thm}\label{cup product cont simple version
thm} Let $(X,g,\omega_g)$ be a \con non\cpt complete K\"ahler
manifold, let $K$ be a closed subset of~$X$, and let $\seq\vphi
j_{j\in J}$ be a countable collection of \cont \plsh \fns on~$X$,
each of which is locally constant on a \nbd of~$K$ in~$X$. Then
for every constant $\epsilon\in(0,1)$, there exists a complete
K\"ahler metric~$h_0$ on~$X$ \st $h_0\geq(1-\epsilon)g$ on~$X$,
$h_0=g$ at each point in~$K$, and for every \con covering space
$\Upsilon\colon\what X\to X$ and every (complete) K\"ahler metric
$h$ on~$\what X$ with $h\geq\hat h_0\equiv\Upsilon^*h_0$, the set
$\what K\equiv\Upsilon\inv(K)$ and the \cont \plsh \fns
$\set{\hat\vphi_j}$ given by $\hat\vphi_j\equiv\Upsilon^*\vphi_j$ for
each $j\in J$ have the following properties:
\begin{enumerate}

\item[(i)] If $\theta$ is a \cinf $(1,0)$-form on~$\what X$ for
which $d\theta\restrict{\what X\sm\what K}\equiv 0$ and
$\theta\restrict{\what X\sm\what K}$ is in $L^2$ with respect
to~$h$, then \(\partial\hat\vphi_j\wedge\theta\equiv 0\) as a
current on~$\what X$ for each $j\in J$. In particular, for each
\concompns~$U$ of $\what X\sm\what K$ on which $\theta$ is not
everywhere zero, $\hat\vphi_j$ is constant on each leaf of the
\holo foliation determined by~$\theta$ in~$U$ for each $j\in J$.

\item[(ii)] If $\theta_1$ and $\theta_2$ are two \cinf
$(1,0)$-forms on~$\what X$ for which $d\theta_k\restrict{\what
X\sm\what K}\equiv 0$ and $\theta_k\restrict{\what X\sm\what K}$
is in $L^2$ with respect to~$h$ for $k=1,2$, and $U$ is a \concomp
of $\what X\sm\what K$ on which $\hat\vphi_j$ is nonconstant for
some $j\in J$, then $\theta_1\wedge\theta_2\equiv 0$ on~$U$.
\end{enumerate}
\end{thm}

Theorem~\ref{cup product cont simple version thm} is a direct
consequence of the following more general version:
\begin{thm}\label{cup product cont thm} Let
$(X,g_0,\omega_{g_0})$ be a \con non\cpt complete Hermitian
manifold of dimension~$n$, let $K$ be a closed subset of~$X$, let
$\seq\vphi j_{j\in J}$ be a countable collection of \cont \plsh
\fns on~$X$, each of which is locally constant on a \nbd of~$K$
in~$X$, and let $\delta\colon X\to(0,1)$ be a \cont \fnns. Then
there exists a nonnegative \cinf \fnns~$\psi$ on~$X$ \st
$\lev\psi\geq-\delta g_0$ on~$X$, $\psi$ is constant on each
\concomp of~$K$, the derivatives of~$\psi$ of all orders vanish at
each point in~$K$, and for every constant $\epsilon\in(0,1)$,
every \con covering space $\Upsilon\colon\what X\to X$, and every
(complete) Hermitian metric $g$ on~$\what X$ with $g\geq\hat
g_0\equiv\Upsilon^*g_0$, the \cinf
\fnns~$\hat\psi\equiv\Upsilon^*\psi$, the (complete) Hermitian
metric $h\equiv g+\epsilon{\cal L}(\hat\psi)$ with associated
$(1,1)$-form~$\omega_h$, the set $\what K\equiv\Upsilon\inv(K)$,
and the \cont \plsh \fns $\set{\hat\vphi_j}$ given by
$\hat\vphi_j\equiv\Upsilon^*\vphi_j$ for each $j\in J$ have the
following properties:
\begin{enumerate}

\item[(i)] If $\beta$ is a \cinf $(1,1)$-form on~$\what X$ for
which $\beta\restrict{\what X\sm\what K}\geq 0$,
$d(\beta\wedge\omega_h^{n-2})\restrict{\what X\sm\what K}\equiv
0$, and $\beta\restrict{\what X\sm\what K}$ is in $L^1$ with
respect to~$h$, then \((\partial\dbar\hat\vphi_j)\wedge\beta\equiv
0\) and \((\partial\hat\vphi_j)\wedge\beta\equiv 0\) as currents
on~$X$ for each $j\in J$. Moreover, if $U\subset\what X$ is a
domain on which $\beta=i\theta\wedge\bar\theta$ for some \cinf
$(1,0)$-form~$\theta$ on~$U$, then
$(\partial\hat\vphi_j)\wedge\theta\equiv 0$ as a current on~$U$
for each $j\in J$. In particular, if this form $\theta$ is a
nontrivial closed \holo $1$-form, then $\hat\vphi_j$ is constant
on each leaf of the \holo foliation determined by~$\theta$ in~$U$
for each $j\in J$.

\item[(ii)] If $\beta_1$ and $\beta_2$ are two \cinf $(1,1)$-forms
on~$\what X$ for which $\beta_k\restrict{\what X\sm\what K}\geq
0$, $d(\beta_k\wedge\omega_h^{n-2})\restrict{\what X\sm\what
K}\equiv 0$, and $\beta_k\restrict{\what X\sm\what K}$ is in $L^1$
with respect to~$h$ for $k=1,2$, $\theta_1$ and $\theta_2$ are
closed \holo $1$-forms on a domain $U\subset\what X$, $\beta_k\restrict
U=i\theta_k\wedge\bar\theta_k$ for $k=1,2$, and
$\hat\vphi_j\restrict U$ is nonconstant for some $j\in J$, then
$\theta_1\wedge\theta_2\equiv 0$ on~$U$.
\end{enumerate}
\end{thm}

We first consider some elementary observations.
\begin{lem}\label{add term to metric bounds it lem}
Let $(X,g,\omega)$ be a Hermitian manifold of dimension~$n$, and
let $\theta$ be a \cinf form of type~$(1,0)$ on~$X$.
\begin{enumerate}

\item[(a)] If $h$ is the Hermitian metric with associated
$(1,1)$-form $\omega_h\equiv\omega+i\,\theta\wedge\bar\theta$,
then \(|\theta|^2_h=(1+|\theta|^2_g)\inv|\theta|^2_g\leq 1\),
\(dV_h=(1+|\theta|^2_g)dV_g\), and
\(|\theta|^2_h\,dV_h=|\theta|^2_g\,dV_g\). Moreover,
\(|v|_h^2=|v|^2_g+2|\theta(v)|^2=|v|^2_g+2|(\real\theta)(v)|^2+2|(\imag\theta)(v)|^2\)
for every real tangent vector $v\in TX$.

\item[(b)] We have \(\|\theta\|^2_{L^2(X,g)}=\int_X
\frac{\sqrt{-1}}{(n-1)!}\,\theta\wedge\bar\theta\wedge\omega^{n-1}\).
Consequently, if $g'$ is any Hermitian metric with associated
$(1,1)$-form~$\omega'$ and
$\theta\wedge\omega=\theta\wedge\omega'$, then
$\|\theta\|_{L^2(X,g)}=\|\theta\|_{L^2(X,g')}$.

\end{enumerate}
\end{lem}
\begin{pf}
For part~(a), observe that $h=g$ at any point $p\in X$ at which
$\theta_p=0$, while at any point $p\in X$ at which $\theta_p\neq
0$, one may verify the first group of equalities by writing $g$
and $h$ in terms of a $g$-orthonormal basis $e_1,\dots,e_n$ for
$T^{1,0}_pX$ with dual basis
\(e_1^*=|\theta_p|\inv_g\theta_p,e_2^*,\dots , e_n^*\). The
verification of the last equality in part~(a) and the verification
of part~(b) are also straightforward. 
\end{pf}

\begin{lem}[cf.~Lemma~5.1 of \cite{NR L2 Castelnuovo}]\label{positive forms linear alg lemma}
Let $\left(\cal V,J\right)$ be a \cpx vector space of dimension
$n>1$, let~$g$ be a Hermitian inner product on $\cal V$ with
associated real skew-symmetric $(1,1)$-form~$\omega$, let
$\alpha$~and~$\beta$ be real skew-symmetric forms of type~$(1,1)$
on $\cal V$ with $\alpha\geq 0$ and $\beta\geq 0$, and let $\eta$
and $\theta$ be $(1,0)$-forms on~$\cal V$. Then
\begin{enumerate}

\item[(i)] $\frac 1{n!}|\beta|_g\omega^n\leq\frac
1{(n-1)!}\beta\wedge\omega^{n-1}\leq\frac {\sqrt
n}{n!}|\beta|_g\omega^n$;

\item[(ii)] $\frac 1{n!}|\alpha\wedge\beta|_g\omega^n\leq \frac
1{(n-2)!}\alpha\wedge\beta\wedge\omega^{n-2} \leq\frac
{\sqrt{n(n-1)/2}}{n!}|\alpha\wedge\beta|_g\omega^n$;

\item[(iii)] $\frac 1{n!}|\eta\wedge\beta|_g^2\omega^n\leq \frac
1{(n-2)!}|\beta|_g
i\eta\wedge\bar\eta\wedge\beta\wedge\omega^{n-2}\leq
\frac{n-1}{n!}|\eta|^2_g|\beta|^2_g\omega^n$; and

\item[(iv)] $\frac 1{n!}|\eta\wedge\theta|_g^2\omega^n =\frac
1{(n-2)!}(i\eta\wedge\bar\eta)\wedge(i\theta\wedge\bar\theta)\wedge\omega^{n-2}
\leq\frac 1{n!}|\eta|^2_g|\theta|^2_g\omega^n$.

\end{enumerate}
\end{lem}
\begin{pf}
We may choose a basis $\zeta_1,\dots,\zeta_n$ for $\left(\cal
V^*\right)^{1,0}$ so that
\[
\omega=\sum i\zeta_j\wedge\bar\zeta_j\qquad\text{and}\qquad
\beta=\sum \lambda_ji\zeta_j\wedge\bar\zeta_j,
\]
where $0\leq\lambda_1\leq\lambda_2\leq\cdots\leq\lambda_n$. We
then have
\[
\alpha=\sum
A_{jk}i\zeta_j\wedge\bar\zeta_k\qquad\text{and}\qquad\eta=\sum
r_j\zeta_j,
\]
where \(A_{jk}=\overline{A_{kj}}\) for all $j,k$.
Thus
$|\beta|_g=\left(\sum\lambda_j^2\right)^{1/2}$ and
\begin{align*}
\beta\wedge\omega^{n-1}&=i^n(n-1)!\sum_{j=1}^n\sum_{k=1}^n
\lambda_j\zeta_j\wedge\bar\zeta_j\wedge\zeta_1\wedge\bar\zeta_1\wedge
\cdots\wedge\what{\zeta_k\wedge\bar\zeta_k}\wedge\cdots
\wedge\zeta_n\wedge\bar\zeta_n\\
&=i^n(n-1)!\sum_{j=1}^n
\lambda_j\zeta_1\wedge\bar\zeta_1\wedge\cdots
\wedge\zeta_n\wedge\bar\zeta_n =\frac
1n\sum_{j=1}^n\lambda_j\,\omega^n,
\end{align*}
and the claim (i) follows.

We also have
\begin{align*}
\alpha\wedge\beta&=i^2\sum_{j\neq k,l}A_{kl}\lambda_j\zeta_k\wedge\bar\zeta_l\wedge\zeta_j\wedge\bar\zeta_j\\
&=i^2\sum_{\substack{{j,k,l\text{ are}}\\{\text{distinct}}}}\lambda_jA_{kl}\zeta_j\wedge\bar\zeta_j\wedge\zeta_k\wedge\bar\zeta_l
+i^2\sum_{j<k}(\lambda_jA_{kk}+\lambda_kA_{jj})\zeta_j\wedge\bar\zeta_j\wedge\zeta_k\wedge\bar\zeta_k.
\end{align*}
Hence
\begin{align*}
\alpha\wedge&\beta\wedge\omega^{n-2}\\
&=i^{n-2}(n-2)!\sum_{j<k}\alpha\wedge\beta\wedge\zeta_1\wedge\bar\zeta_1\wedge\cdots\wedge
\what{\zeta_j\wedge\bar\zeta_j}\wedge\cdots\wedge
\what{\zeta_k\wedge\bar\zeta_k}
\wedge\cdots\wedge\zeta_n\wedge\bar\zeta_n\\
&=i^n(n-2)!\sum_{j<k}(\lambda_jA_{kk}+\lambda_kA_{jj})
\zeta_1\wedge\bar\zeta_1\wedge\cdots\wedge\zeta_n\wedge\bar\zeta_n\\
&=\frac 1{n(n-1)}\left(\sum_{j<k}(\lambda_jA_{kk}+\lambda_kA_{jj})\right)\omega^n,
\end{align*}
and
\begin{align*}
|\alpha\wedge\beta|^2_g&=\sum_{\substack{{j,k,l\text{ are}}\\{\text{distinct}}}}\lambda_j^2|A_{kl}|^2
+\sum_{j<k}(\lambda_jA_{kk}+\lambda_kA_{jj})^2\leq\sum_{\substack{{j,k,l\text{ are}}\\{\text{distinct}}}}\lambda_j^2A_{kk}A_{ll}
+\sum_{j<k}(\lambda_jA_{kk}+\lambda_kA_{jj})^2\\
&\leq\left(\sum_{j<k}(\lambda_jA_{kk}+\lambda_kA_{jj})\right)^2\leq\frac{n(n-1)}2\sum_{j<k}(\lambda_jA_{kk}+\lambda_kA_{jj})^2
\leq\frac{n(n-1)}2|\alpha\wedge\beta|^2_g,
\end{align*}
so (ii) follows.

From the above, we have
\[
i\eta\wedge\bar\eta\wedge\beta\wedge\omega^{n-2}= \frac
1{n(n-1)}\left(\sum_{j\neq k}\lambda_j|r_k|^2\right)\omega^n,
\]
while
\begin{align*}
|\eta\wedge\beta|_g^2&= \left|\sum_{j\neq k}
\lambda_jr_k\zeta_k\wedge\zeta_j\wedge\bar\zeta_j\right|^2_g
=\sum_{j\neq k} \lambda_j^2|r_k|^2\leq |\beta|_g\sum_{j\neq k}
\lambda_j|r_k|^2\\& \leq|\beta|^2_g\sum_{j\neq k}|r_k|^2
=(n-1)|\eta|^2_g|\beta|^2_g,
\end{align*}
so (iii) follows.

Finally, letting $\beta=i\theta\wedge\bar\theta$, we get
$\lambda_j=0$ for $j=1,\dots,n-1$, $\lambda_n=|\theta|^2_g$,
\[
(i\eta\wedge\bar\eta)\wedge(i\theta\wedge\bar\theta)\wedge\omega^{n-2}
=\frac
1{n(n-1)}\left(\sum_{k=1}^{n-1}|\theta|^2_g|r_k|^2\right)\omega^n,
\]
and
\begin{align*}
|\eta\wedge\theta|_g^2&= \left|\sum_{k=1}^{n-1}
r_k|\theta|_g\zeta_k\wedge\zeta_n\right|^2_g =\sum_{k=1}^{n-1}
|\theta|^2_g|r_k|^2\leq|\eta|^2_g|\theta|^2_g,
\end{align*}
so (iv) follows.

\end{pf}

Given a Hermitian inner product~$g$ with
associated real skew-symmetric $(1,1)$-form~$\omega$ on a \cpx vector space $\left(\cal
V,J\right)$ of dimension~$n$, and a real skew-symmetric
form~$\alpha$ of type~$(1,1)$ on~$\cal V$, we have the orthogonal decomposition
$\cal V^{1,0}=\cal V^{1,0}_-\oplus\cal V_0\oplus\cal V^{1,0}_+$,
where $\cal V^{1,0}_-$ ($\cal V_0$, $\cal V^{1,0}_+$) is the sum
of the eigenspaces for the positive eigenvalues (respectively, the
eigenspace for the zero eigenvalue, the sum of the eigenspaces for
the negative eigenvalues). Letting \(\text{pr}_{+}\colon\cal
V^{1,0}\to\cal V^{1,0}_+\) and \(\text{pr}_{-}\colon\cal
V^{1,0}\to\cal V^{1,0}_-\) be the corresponding orthogonal
projections, we get $\alpha=\alpha^+-\alpha^-$, where $\alpha^+$
and $\alpha^-$, the \emph{positive part of~$\alpha$} and
\emph{negative part of~$\alpha$}, respectively, are the
nonnegative $(1,1)$-forms given by
\[
\alpha^+(u,\bar v)=\alpha(\text{pr}_{+}u,
\overline{\text{pr}_{+}v})\qquad\text{and}\qquad\alpha^{-}(u,\bar
v)=-\alpha(\text{pr}_{-}u,
\overline{\text{pr}_{-}v})\qquad\forall\,u,v\in\cal V^{1,0}.
\]
If $\seq\zeta j$ is a basis for $\left(\cal V^*\right)^{1,0}$ in
which \(\omega=\sum i\zeta_j\wedge\overline{\zeta_j}\) and
\(\alpha=\sum\lambda_ji\zeta_j\wedge\overline{\zeta_j}\), then
\[
\alpha^+=\sum\lambda^+_ji\zeta_j\wedge\overline{\zeta_j}
\qquad\text{and}\qquad
\alpha^-=\sum\lambda^-_ji\zeta_j\wedge\overline{\zeta_j}.
\]
\begin{lem}\label{Pos neg parts of cont forms are cont lem}
Let $\alpha$ be a real $(1,1)$-form on a Hermitian manifold
$(X,g,\omega)$. If $\alpha$ is \cont (measurable), then the
associated positive and negative parts, $\alpha^+$ and $\alpha^-$,
are \cont (respectively, measurable).
\end{lem}
\begin{pf}
Let $n=\dim X$. If $\alpha$ is \cont but $\alpha^+$ is not, then
there exist a point $p\in X$, tangent vectors
$u,v\in\left(T_pX\right)^{1,0}$, sequences $\seq x\nu$ in~$X$ and
$\seq u\nu$ and $\seq v\nu$ in $\left(TX\right)^{1,0}$, and a
positive constant~$\epsilon$ \st
$u_\nu,v_\nu\in\left(T_{x_\nu}X\right)^{1,0}$ for each~$\nu$,
$x_\nu\to p$, $u_\nu\to u$, and $v_\nu\to v$, but
\[
|\alpha^+(u_\nu,\bar v_\nu)-\alpha^+(u,\bar
v)|\geq\epsilon\qquad\forall\,\nu=1,2,3,\dots.
\]
For each~$\nu$, we may fix an (orthonormal) basis
$\set{\zeta\ssp\nu_j}_{j=1}^n$ for $(T_{x_\nu}^*X)^{1,0}$ \st
\[
\omega_{x_\nu}=\sum_{j=1}^ni\zeta\ssp\nu_j\wedge
\overline{\zeta\ssp\nu_j}\qquad\text{and}\qquad
\alpha_{x_\nu}=\sum_{j=1}^n\lambda\ssp\nu_ji\zeta\ssp\nu_j\wedge
\overline{\zeta\ssp\nu_j},
\]
where $\lambda\ssp\nu_1\leq\cdots\leq\lambda\ssp\nu_n$. By the
continuity of~$\alpha$, the eigenvalues $\set{\lambda\ssp\nu_j}$
are uniformly bounded, so by replacing the above with a suitable
subsequence, we may assume that
\[
\zeta\ssp\nu_j\to\zeta_j\qquad\text{and}\qquad\lambda\ssp\nu_j\to\lambda_j
\qquad\forall\, j=1,\dots,n,
\]
and hence that \(\omega_p=\sum i\zeta_j\wedge\bar\zeta_j\) and
\(\alpha_p=\sum\lambda_ji\zeta_j\wedge\bar\zeta_j\). However, this
then implies that
\[
\alpha^+(u_\nu,\bar
v_\nu)=\sum(\lambda_j\ssp\nu)^+i\zeta\ssp\nu_j(u_\nu)
\overline{\zeta\ssp\nu_j(v_\nu)}\longrightarrow
\sum\lambda_j^+i\zeta_j(u)\overline{\zeta_j(v)}=\alpha^+(u,\bar
v).
\]
Thus we have arrived at a contradiction, and hence continuity
of~$\alpha$ implies that of~$\alpha^+$.

For $\alpha$ a measurable form, there exists a sequence of \cont
real $(1,1)$-forms $\seq\alpha\nu$ converging to~$\alpha$ almost
everywhere in~$X$, and an argument similar to the above shows that
$\alpha^+_\nu\to\alpha^+$~a.e. in~$X$. Hence $\alpha^+$ is also
measurable.
\end{pf}

For the proof of Theorem~\ref{cup product cont thm}, after forming \cinf approximations of the
given \plsh \fns and modifying
the metric, we will apply the following:
\begin{lem}\label{cinf approx cup product pass to limit lem}
Let $\seq\alpha\nu_{\nu=1}^\infty$ and $\beta$ be \cont real
$(1,1)$-forms on a Hermitian manifold $(X,g,\omega)$ of
dimension~$n$ \st $\beta\geq 0$ and
$\alpha_\nu\wedge\beta\wedge\omega^{n-2}\to 0$ in~$L^1_{\text{\rm
loc}}$.
\begin{enumerate}

\item[(a)] If $\alpha_\nu^-\wedge\beta\to 0$ in $L^1_{\text{\rm
loc}}$, then $\alpha_\nu\wedge\beta\to 0$ in $L^1_{\text{\rm
loc}}$.

\item[(b)] If for each~$\nu$,
$\alpha_\nu=i\eta_\nu\wedge\bar\eta_\nu$ for some \cont
$(1,0)$-form~$\eta_\nu$, then $\eta_\nu\wedge\beta\to 0$ in
$L^2_{\text{\rm loc}}$.

\item[(c)] If $\beta=i\theta\wedge\bar\theta$ for some \cont
$(1,0)$-form~$\theta$ and
$\alpha_\nu=i\eta_\nu\wedge\bar\eta_\nu$ for some \cont
$(1,0)$-form~$\eta_\nu$  for each~$\nu$, then $\eta_\nu\wedge\theta\to 0$ in
$L^2_{\text{\rm loc}}$.

\end{enumerate}
\end{lem}
\begin{pf}
Assuming that $\alpha_\nu^-\wedge\beta\to 0$ in $L^1_{\text{\rm
loc}}$, Lemma~\ref{positive forms linear alg lemma} implies that
$\alpha_\nu^-\wedge\beta\wedge\omega^{n-2}\to 0$ in
$L^1_{\text{\rm loc}}$, and hence that
\[
\alpha_\nu^+\wedge\beta\wedge\omega^{n-2}=
\alpha_\nu\wedge\beta\wedge\omega^{n-2}
+\alpha_\nu^-\wedge\beta\wedge\omega^{n-2}\to 0\qquad\text{in
}L^1_{\text{\rm loc}}.
\]
Applying the lemma once more, we get $\alpha_\nu^+\wedge\beta\to
0$ in $L^1_{\text{\rm loc}}$, and part~(a) follows. Parts~(b)
and~(c) follow immediately from parts~(iii) and~(iv),
respectively, of Lemma~\ref{positive forms linear alg lemma}.
\end{pf}

After we obtain the main parts of Theorem~\ref{cup product cont thm}, the following
elementary observations will give the remaining conclusions:
\begin{lem}\label{plsh fn holo 1 form wedge 0 lem}
Let $\vphi$ be a nonconstant real-valued \cont \fn on a \con \cpx
manifold~$X$.
\begin{enumerate}
\item[(a)] For any nontrivial closed \holo $1$-form~$\theta$
on~$X$, the following are equivalent:
\begin{enumerate}

\item[(i)] As a current, $\partial\vphi\wedge\theta\equiv 0$.

\item[(ii)] As a current,
$\partial\vphi\wedge\theta\wedge\bar\theta\equiv 0$.

\item[(iii)] The restriction of $\vphi$ to each leaf of the \holo
foliation determined by~$\theta$ is constant  (equivalently,
for every \holo \fnns~$f$ with $\theta=df$ on an open set~$U$,
$\vphi$ is constant on each level of~$f$).

\end{enumerate}

\item[(b)] If $\theta_1$ and $\theta_2$ are closed \holo $1$-forms
on $X$ and $\partial\vphi\wedge\theta_j\equiv 0$ for $j=1,2$, then
$\theta_1\wedge\theta_2\equiv 0$.

\end{enumerate}
\end{lem}
\begin{pf}
Let $n\equiv\dim X$. For the proof of (a), given a point $p\in X$
at which $\theta_p\neq 0$, we may fix a local \holo coordinate
\nbd $(U,(z_1,\dots,z_n))$ of~$p$ in which $U$ is a coordinate
polydisk and $\theta\restrict U=dz_1$. For each $j=2,\dots,n$, and
for each \fn $u\in\cal D(U)$, we then have
\begin{align*}
\langle\partial\vphi\wedge\theta\wedge\bar\theta,&u\,dz_2\wedge
d\bar z_2\wedge\cdots\wedge\what{dz_j}\wedge d\bar
z_j\wedge\cdots\wedge dz_n\wedge d\bar z_n\rangle
\\
&=\langle\partial\vphi\wedge\theta,u\,d\bar z_1\wedge dz_2\wedge
d\bar z_2\wedge\cdots\wedge\what{dz_j}\wedge d\bar
z_j\wedge\cdots\wedge dz_n\wedge d\bar z_n\rangle
\\
&=\langle\partial\vphi,u\,dz_1\wedge d\bar
z_1\wedge\cdots\wedge\what{dz_j}\wedge d\bar z_j\wedge\cdots\wedge
dz_n\wedge d\bar z_n\rangle\\
&=-\int_U\vphi\pdof u{z_j}\,dz_1\wedge d\bar z_1\wedge\cdots\wedge
dz_n\wedge d\bar z_n.
\end{align*}
Thus the conditions $\partial\vphi\wedge\theta\equiv 0$ and
$\partial\vphi\wedge\theta\wedge\bar\theta\equiv 0$ in~$U$ are
each equivalent to the condition
\[
\left(\pdof\vphi{\bar z_j}\right)_{\text{distr.}}
=\left(\pdof\vphi{z_j}\right)_{\text{distr.}}=0 \qquad\text{for
}j=2,\dots,n.
\]
Hence the above conditions are equivalent to the condition that
$\vphi\restrict U$ is a \fn of~$z_1$, i.e., that $\vphi$ is
constant on the leaves of the \holo foliation determined by
$\theta\restrict U$.

Given an arbitrary point $p\in X$, we may choose a \con \rel \cpt
\nbdns~$U$ and a bounded \holo \fn $f$ on~$U$ \st $\theta\restrict
U=df$, $f(p)=0$, the fiber $L\equiv f\inv(0)$ through~$p$ is
\conns, and $\theta$ is nonvanishing on $U\sm L$.
By continuity of intersections (see \cite{Stein} or
\cite{Tworzewski-Winiarski Cont of intersect} or Theorem~4.23 in
\cite{ABCKT}), after fixing a sequence of points $\seq x\nu$ in
$U\sm L$ converging to~$p$, letting $L_\nu$ be the level of~$f$
through~$x_\nu$ for each~$\nu$, and passing to a subsequence, we
get a sequence of \anal sets $\seq L\nu$ converging to~$L$ in~$U$.
By the above, if $\partial\vphi\wedge\theta\equiv 0$ or
$\partial\vphi\wedge\theta\wedge\bar\theta\equiv 0$ in~$U$, then
$\vphi\restrict{L_\nu}$ is constant for each~$\nu$, and hence
$\vphi\restrict L$ is constant. Conversely, assuming that $\vphi$
is constant on each level of~$f$, let us fix a constant
$m<\inf_U\vphi$, and let us set
\[
\vphi_\nu\equiv\max(\vphi+\nu\inv\log|f|,m)\qquad\text{for
}\nu=1,2,3,\dots.
\]
Then $\seq\vphi\nu$ is a uniformly bounded sequence of \cont \fns
converging pointwise almost everywhere to $\vphi\restrict U$.
Moreover, each term is constant on a \nbd of~$L$ as well as on
each level of~$f$. Applying the previous observations and the
dominated convergence theorem, we see that for every form
$\alpha\in\cal D^{n-2,n}(U)$,
\[
0=\langle\partial\vphi_\nu\wedge\theta,\alpha\rangle\to
\langle\partial\vphi\wedge\theta,\alpha\rangle.
\]
Thus $\partial\vphi\wedge\theta\equiv 0$, and hence
$\partial\vphi\wedge\theta\wedge\bar\theta\equiv 0$, on~$U$. The claim~(a)
now follows.

For the proof of (b), observe that each point in the complement of
the zero set~$Z$ of $\theta_1\wedge\theta_2$ admits a local \holo
coordinate polydisk \nbd $(U,(z_1,\dots,z_n))$ \st
$dz_1=\theta_1\restrict U$ and $dz_2=\theta_2\restrict U$. The
restriction $\vphi\restrict U$ is then both a \fn of $z_1$ and a
\fn of $z_2$, and is therefore constant. Hence the nonconstant
\fnns~$\vphi$ is constant on the \con set $X\sm Z$, which is
either dense or empty, so we must have $Z=X$.
\end{pf}

The main step in the proof of Theorem~\ref{cup product cont thm} is the
following:
\begin{lem}\label{cup product preliminary lem}
Let $(X,g,\omega)$ be a \con non\cpt complete Hermitian manifold
of dimension~$n$, let $K$ be a closed subset of~$X$, let $\rho$ be
a real-valued \cinf \fn on~$X$, let $\beta$ be a \cinf
$(1,1)$-form on~$X$, let $\gamma\equiv i\partial\dbar\rho\wedge\beta\wedge\omega^{n-2}$,
and let $R$ be a bounded positive
\cont \fn on~$X$. Assume that
\begin{enumerate}

\item[(i)] $|d\rho|_g$ is bounded on~$X$,
$i\partial\dbar\rho\geq-R\omega$ on $X$, and $\rho$ is locally constant on some \nbd
of~$K$; and

\item[(ii)] We have $\beta\restrict{X\sm K}\geq 0$,
$d(\beta\wedge\omega^{n-2})\restrict{X\sm K}\equiv 0$, and
$\beta\restrict{X\sm K}$ is in $L^1$.

\end{enumerate}
Then
\[
\int_X\gamma^+
=\int_X\gamma^-
\leq\int_{X\sm K}\frac R{\sqrt n}|\beta|_g\,\omega^n<\infty
\]
(where at each point, $\gamma^\pm\equiv(\gamma/\omega_0^n)^\pm\cdot\omega_0^n$
for any positive $(1,1)$-form~$\omega_0$).
In particular, $\gamma$ is integrable, $\int_X\gamma=0$, and
\[
\|\gamma\|_{L^1}=\int_X\gamma^++\int_X\gamma^-\leq\int_{X\sm K}\frac {2R}{\sqrt n}|\beta|_g\,\omega^n.
\]
\end{lem}
\begin{pf}
As in~\cite{Gaffney}, fixing a point $p\in X$ and setting
\[
\tau(s)\equiv\left\{
\begin{aligned}
1&\quad\text{if }s\leq 1\\
2-s&\quad\text{if }1<s<2\\
0&\quad\text{if }2\leq s
\end{aligned}\right.
\]
and
\[
\tau_r(x)\equiv\tau\left(\frac{\text{dist}(p,x)}{r}\right)\qquad\forall\,x\in
X,r>0,
\]
we get a collection of nonnegative Lipschitz \cont \fns $\seq\tau
r_{r>0}$ \st for each $r>0$, we have $0\leq\tau_r\leq 1$ on $X$,
$\tau_r\equiv 1$ on $B(p;r)$, $\tau_r\equiv 0$ on~$X\sm B(p;2r)$,
and $|d\tau_r|_g\leq 1/r$. We then have
\[
\int_X2\tau_r\gamma=\int_X\tau_rdd^c\rho
\wedge\beta\wedge\omega^{n-2}\\
=-\int_{(B(p;2r)\sm B(p;r))\sm K}d\tau_r\wedge
d^c\rho\wedge\beta\wedge\omega^{n-2},
\]
where $d^c=-i(\partial-\dbar)$. For some positive constant $C$, we
have
\[
\left|d\tau_r\wedge d^c\rho\wedge\beta\wedge\omega^{n-2}\right|_g
\leq\frac Cr|\beta|_g,
\]
and hence
\(\int_X\tau_r\gamma\to 0\) as $r\to\infty$.
Since $i\partial\dbar\rho\geq-R\omega$, Lemma~\ref{positive forms
linear alg lemma} implies that on $X\sm K$,
\[
\gamma
\geq-R\,\beta\wedge\omega^{n-1}\geq-\frac{R}{\sqrt
n}|\beta|_g\,\omega^n.
\]
Applying the monotone convergence theorem, we
get
\[
\int_X\gamma^+=
\lim_{r\to\infty}\int_X\left(\tau_r\gamma\right)^+
=\lim_{r\to\infty}\int_X\left(\tau_r\gamma\right)^-
=\int_X\gamma^-
\leq\int_{X\sm K}\frac R{\sqrt n}|\beta|_g\,\omega^n.
\]
\end{pf}

\begin{pf*}{Proof of Theorem~\ref{cup product cont thm}}
Let $\seq\delta{\nu j}_{\nu\in\N,j\in J}$ be a family of bounded positive
\cont \fns on~$X$, the elements of which we will later choose to
be sufficiently small. For each $\nu\in\N$ and each $j\in J$, the
\cinf approximation theorem of Greene and Wu (Theorem~\ref{cinf
approximation thm}) provides a \cinf \fn $\vphi_{\nu j}$ \st
$\vphi_{\nu j}=\vphi_j$ on a \nbd of~$K$, and
\[
\vphi_j\leq\vphi_{\nu j}<\vphi_j
+\delta_{\nu j}
\qquad\text{and}\qquad\lev{\vphi_{\nu j}}\geq-\delta_{\nu j}g_0
\qquad\text{on }X.
\]
So that we may work with positive \fnsns, let us set 
\(\psi_{\nu j}\equiv e^{\vphi_{\nu j}}\) 
for each $\nu\in\N$ and $j\in J$. We
will obtain the required \cinf \fnns~$\psi$ as a weighted sum of
the squares of the \fns $\seq\psi{\nu j}$. By Lemma~\ref{add term
to metric bounds it lem}, the liftings of the \fns $\seq\psi{\nu
j}$ to a covering space will then have bounded gradient with
respect to the associated modified metric. While these \fns are
not \plshns, Lemma~\ref{cup product preliminary lem}
will imply that the cup product property holds in an approximate
sense, and we will then pass to a limit.

For each $\nu\in\N$ and $j\in J$, we have
\[
i\partial\dbar\psi_{\nu j}
=e^{\vphi_{\nu j}}\left(i\partial\dbar\vphi_{\nu j}
+i\partial\vphi_{\nu j}\wedge\dbar\vphi_{\nu j}\right)
\geq -e^{\vphi_j+\delta_{\nu j}}\delta_{\nu j}\omega_{g_0}
\]
and
\begin{align*}
i\partial\dbar\psi_{\nu j}^2 
&=2\psi_{\nu j}i\partial\dbar\psi_{\nu j}
+2i\partial\psi_{\nu j}\wedge\dbar\psi_{\nu j}\\
&\geq-2e^{2\vphi_j+2\delta_{\nu j}}
\delta_{\nu j}\omega_{g_0}
+2i\partial\psi_{\nu j}\wedge\dbar\psi_{\nu j}
\geq-2e^{2\vphi_j+2\delta_{\nu j}}\delta_{\nu j}\omega_{g_0}.
\end{align*}
Hence, choosing $\delta_{\nu j}$ so small that 
\(0<2e^{2\vphi_j+2\delta_{\nu j}}\delta_{\nu j}
<\delta<1\),
we get
\[
i\partial\dbar\psi_{\nu j}^2 
\geq-\delta\omega_{g_0}
+2i\partial\psi_{\nu j}\wedge\dbar\psi_{\nu j}\geq-\delta\omega_{g_0}.
\]

If  $\seq\epsilon{\nu j}_{\nu\in\N,j\in J}$ is a family of
sufficiently small positive constants (depending on the choice of
$\seq\delta{\nu j}$), then  $\sum_{\nu,j}\epsilon_{\nu j}<1$;
the series
\[
\sum_{\nu\in\N,j\in J}\epsilon_{\nu j}\psi_{\nu j}^2=
\sum_{\nu\in\N,j\in J}\epsilon_{\nu j}e^{2\vphi_{\nu j}}
\]
converges to a nonnegative \cinf \fnns~$\psi$; for each
$k=0,1,2,\dots$, any term-by-term $k$th order derivative series
for the above series converges uniformly and absolutely on \cpt
subsets of~$X$ to the corresponding $k$th order derivative
of~$\psi$; and each of the series 
\[
\sum_{\nu,j}\epsilon_{\nu j}\psi_{\nu j}i\partial\dbar\psi_{\nu j}
\qquad\text{and}\qquad
\sum_{\nu,j}\epsilon_{\nu j}i\partial\psi_{\nu j}\wedge\dbar\psi_{\nu j}
\]
converges uniformly on \cpt subsets of~$X$ to a \cont real  $(1,1)$-form.
Hence
\[
i\partial\dbar\psi=\sum_{\nu,j}\epsilon_{\nu j}i\partial\dbar\psi_{\nu j}^2
\geq-\delta\omega_{g_0}
+\sum_{\nu,j}2\epsilon_{\nu j}i\partial\psi_{\nu j}\wedge\dbar\psi_{\nu j}
\geq-\delta\omega_{g_0}.
\]

Let us now fix a \con covering space $\Upsilon\colon\what X\to X$,
a Hermitian metric $g\geq\hat g_0\equiv\Upsilon^*g_0$, and a
constant $\epsilon\in(0,1)$. Let $\hat\psi\equiv\Upsilon^*\psi$,
let $h\equiv g+\epsilon\cal L(\hat\psi)$, let $\what
K\equiv\Upsilon\inv(K)$, let $\hat\vphi_j\equiv\Upsilon^*\vphi_j$
for each $j\in J$, 
let $\hat\vphi_{\nu j}\equiv\Upsilon^*\vphi_{\nu j}$, 
$\hat\psi_{\nu j}\equiv\Upsilon^*\psi_{\nu j}=e^{\hat\vphi_{\nu j}}$, and
$\hat\delta_{\nu j}\equiv\Upsilon^*\delta_{\nu j}$ for all
$\nu\in\N$ and $j\in J$, let $\hat\delta\equiv\Upsilon^*\delta$,
and let $\omega_{\hat g_0}$, $\omega_g$, and
$\omega_h=\omega_g+\epsilon\,i\partial\dbar\hat\psi$ be the
$(1,1)$-forms associated to $\hat g_0$, $g$, and $h$,
respectively. Suppose $\beta$ is a \cinf $(1,1)$-form on~$\what X$ for
which $\beta\restrict{\what X\sm\what K}\geq 0$,
$d(\beta\wedge\omega_h^{n-2})\restrict{\what X\sm\what K}\equiv
0$, and $\beta\restrict{\what X\sm\what K}$ is in $L^1$ with
respect to~$h$. For each $\nu\in\N$ and $j\in J$, we have
\[
\omega_h-\epsilon\epsilon_{\nu j}i\partial\hat\psi_{\nu
j}\wedge\dbar\hat\psi_{\nu j}\geq\omega_{\hat
g_0}-\epsilon\hat\delta\omega_{\hat g_0}>0,
\]
and therefore, by Lemma~\ref{add term to metric bounds it lem},
\[
|d\hat\psi_{\nu j}|_h\leq\sqrt{\frac 2{\epsilon\epsilon_{\nu j}}}.
\]
Since
\[
i\partial\dbar\hat\psi_{\nu j}\geq-e^{\hat\vphi_j
+\hat\delta_{\nu j}}\hat\delta_{\nu j}\omega_{\hat
g_0}\geq-e^{\hat\vphi_j+\hat\delta_{\nu j}}
\hat\delta_{\nu j}(1-\epsilon)\inv\omega_{h},
\]
Lemma~\ref{cup product preliminary lem} implies that
\[
\int_{\what X}\left[i\partial\dbar\hat\psi_{\nu
j}\wedge\beta\wedge\omega_h^{n-2}\right]^+
=\int_{\what
X}\left[i\partial\dbar\hat\psi_{\nu
j}\wedge\beta\wedge\omega_h^{n-2}\right]^- \leq\int_{\what
X\sm\what K}\frac{e^{\hat\vphi_j+\hat\delta_{\nu j}}\hat\delta_{\nu
j}}{(1-\epsilon)\sqrt n}|\beta|_h\omega_h^n.
\]
Thus if we choose the sequence of \fns $\seq\delta{\nu
j}_{\nu=1}^\infty$ for each~$j\in J$ so that
$e^{\vphi_j+\delta_{\nu j}}\delta_{\nu j}\to 0$ uniformly on~$X$
as $\nu\to\infty$, then the sequence of \cinf forms
$\set{i\partial\dbar\hat\psi_{\nu
j}\wedge\beta\wedge\omega_h^{n-2}}_{\nu=1}^\infty$ must converge
to~$0$ in~$L^1$ (with respect to any Hermitian metric, since
the forms are of type $(n,n)$). Applying Lemma~\ref{positive
forms linear alg lemma}, we also get
\begin{align*}
\int_{\what X}\frac 1{n!}\left|(i\partial\dbar\hat\psi_{\nu
j})^-\wedge\beta\right|_h\,\omega_h^n&\leq\int_{\what X}\frac
1{(n-2)!}(i\partial\dbar\hat\psi_{\nu
j})^-\wedge\beta\wedge\omega_h^{n-2}\\
&\leq\int_{\what X\sm\what K}\frac {e^{\hat\vphi_j+\hat\delta_{\nu
j}}\hat\delta_{\nu
j}}{(n-2)!(1-\epsilon)}\beta\wedge\omega_h^{n-1}\\
&\leq\int_{\what X\sm\what K}\frac {e^{\hat\vphi_j+\hat\delta_{\nu
j}}\hat\delta_{\nu j}}{(n-2)!(1-\epsilon)\sqrt
n}|\beta|_h\omega_h^n\to 0,
\end{align*}
and hence, by Lemma~\ref{cinf approx cup product pass to limit
lem}, the sequence of \cinf forms
$\set{i\partial\dbar\hat\psi_{\nu j}\wedge\beta}_{\nu=1}^\infty$
converges to~$0$ in~$L^1_{\text{loc}}$. In order to obtain the
same property for~$\set{\hat\vphi_{\nu j}}$, observe that
\[
\hat\psi_{\nu j}\inv i\partial\dbar\hat\psi_{\nu
j}\wedge\beta=\left[\left(i\partial\dbar\hat\vphi_{\nu
j}\right)^++i\partial\hat\vphi_{\nu j}\wedge\dbar\hat\vphi_{\nu
j}\right]\wedge\beta-\left(i\partial\dbar\hat\vphi_{\nu
j}\right)^-\wedge\beta\qquad\forall\,\nu,
\]
and since the \fns $\set{\hat\vphi_{\nu j}}_{\nu=1}^\infty$ are
uniformly bounded on \cpt sets, the left hand side converges
to~$0$ in $L^1_{\text{loc}}$.  Since
$\seq\delta{\nu j}$ converges to~$0$ uniformly on \cpt sets, the sequence
\(\set{\left(i\partial\dbar\hat\vphi_{\nu
j}\right)^-\wedge\beta}_{\nu=1}^\infty\) must also converge to~$0$
in $L^1_{\text{loc}}$, and therefore so must the sequence
\[
\left\{\left[\left(i\partial\dbar\hat\vphi_{\nu
j}\right)^++i\partial\hat\vphi_{\nu j}\wedge\dbar\hat\vphi_{\nu
j}\right]\wedge\beta\right\}.
\]
Hence by Lemma~\ref{positive forms linear alg lemma}, the sequence
\[
\left\{\left[\left(i\partial\dbar\hat\vphi_{\nu
j}\right)^++i\partial\hat\vphi_{\nu j}\wedge\dbar\hat\vphi_{\nu
j}\right]\wedge\beta\wedge\omega_h^{n-2}\right\},
\]
and therefore, the sequences
\[
\left\{\left(i\partial\dbar\hat\vphi_{\nu
j}\right)^+\wedge\beta\wedge\omega_h^{n-2}\right\}
\qquad\text{and}\qquad \left\{(i\partial\hat\vphi_{\nu
j}\wedge\dbar\hat\vphi_{\nu
j})\wedge\beta\wedge\omega_h^{n-2}\right\},
\]
must converge to $0$ in $L^1_{\text{loc}}$. Applying
Lemma~\ref{cinf approx cup product pass to limit lem}, we see that
the sequences
\[
\left\{(\partial\dbar\hat\vphi_{\nu j})\wedge\beta\right\}
\qquad\text{and}\qquad \left\{(\partial\hat\vphi_{\nu
j})\wedge\beta\right\},
\]
must also converge to $0$ in $L^1_{\text{loc}}$ and
$L^2_{\text{loc}}$, respectively. Since $\hat\vphi_{\nu
j}\to\hat\vphi_j$ uniformly on \cpt sets, it follows that
$(\partial\dbar\hat\vphi_j)\wedge\beta=0$ and
$(\partial\hat\vphi_j)\wedge\beta=0$ for each $j\in J$. Moreover,
if on some open set~$U$, $\beta=i\theta\wedge\bar\theta$ for some
\cinf $(1,0)$-form~$\theta$, then Lemma~\ref{cinf approx cup
product pass to limit lem} implies that for each $j\in J$,
$\left\{(\partial\hat\vphi_{\nu j})\wedge\theta\right\}$ converges
to $0$ in $L^2_{\text{loc}}$ in~$U$, and hence
$(\partial\hat\vphi_j)\wedge\theta=0$. The remaining claims of the
theorem follow from Lemma~\ref{plsh fn holo 1 form wedge 0
lem}.\end{pf*}

\section{Strict plurisubharmonicity
and the Bochner--Hartogs property}\label{bh plsh two pos ev sect}

The main goal of this section is a proof of Theorem~\ref{BH plsh fn cont two pos ev thm},
which we will obtain as an application
of Theorem~\ref{cup product cont very simple version intro thm}. 
The two main ingredients of the proof are Theorem~\ref{cup product cont very simple version intro thm} and the following:
\begin{thm}[see Theorem~3.3 of \cite{NR BH Bounded Geometry}]\label{Bochner Hartogs no L2 forms thm}
Let $X$ be a \con non\cpt hyperbolic complete K\"ahler manifold
with no nontrivial $L^2$ \holo $1$-forms.
\begin{enumerate}

\item[(a)] For every compactly supported $\dbar$-closed \cinf
form~$\alpha$ of type~$(0,1)$ on~$X$, there exists a bounded \cinf
\fnns~$\beta$ with finite energy on~$X$ \st $\dbar\beta=\alpha$
on~$X$ and $\beta$ vanishes on every hyperbolic end $E$ of~$X$
that is contained in $X\sm\supp\alpha$.

\item[(b)] In any ends decomposition $X\sm K=E_1\cup\cdots\cup
E_m$, exactly one of the ends, say~$E_1$, is hyperbolic, and
moreover, every \holo \fn on~$E_1$ admits a (unique) extension to
a \holo \fn on~$X$.

\item[(c)] If $e(X)=1$ (equivalently, every end of~$X$ is
hyperbolic), then $H^1_c(X,\ol)=0$.

\end{enumerate}
\end{thm}
For the proof of Theorem~\ref{BH plsh fn cont two pos ev thm}, we
will apply Theorem~\ref{Bochner Hartogs no L2 forms thm} to
suitable sublevels of~$\vphi$. In order to complete the K\"ahler
metric on these sublevels, we will apply Lemma~\ref{completing
metric on sublevel lem} and Lemma~\ref{completing metric on
sublevel ends lem} below, which are contained implicitly in the
work of Nakano~\cite{Nk}, Greene and Wu~\cite{Greene-Wu}, and
Demailly~\cite{Dem}:
\begin{lem}[cf.~\cite{Nk} and \cite{Dem}]\label{completing metric on sublevel lem}
Let $(X,g,\omega_g)$ be a \con Hermitian manifold with distance
\fn $d_g(\cdot,\cdot)$, let $\epsilon\in(0,1)$, let $\delta\colon
X\to(0,\epsilon]$ be a \cont \fnns, and let $U$ be an open
subset.
\begin{enumerate}

\item[(a)] If $Y$ is a domain in~$X$ for which the restricted
distance \fn $d_g\restrict{\overline Y\sm U}$ is complete and
$\psi$ is a \cinf \fn on~$Y$ \st $\psi<0$, $\psi\to 0$
at~$\partial Y$, and \(\lev\psi\geq\delta\,\psi g\), then
\[
h\equiv g+\lev{-\log(-\psi)}\geq(1-\delta)g
\]
is a Hermitian metric on~$Y$ (which is K\"ahler if~$g\restrict Y$
is K\"ahler) for which the restriction $d_h\restrict{Y\sm U}$ of
the associated distance \fn $d_h(\cdot,\cdot)$ in~$Y$ is complete.

\item[(b)] Suppose $\psi$ is a positive \cinf \fn on~$X$ 
\st $2\psi\lev\psi\geq -\delta g$ on~$X$ and \st there
exists a nonempty closed \con set $K\subset X$ for which
$\psi$ is bounded on~$K$ and for every $a>\sup_K\psi$, the
restriction $d_g\restrict{\overline Y\sm U}$
of~$d_g(\cdot,\cdot)$, where $Y$ is the \concomp of $\setof{x\in
X}{\psi(x)<a}$ containing~$K$, is complete. Then
\[
h\equiv g+\lev{\psi^2}\geq(1-\delta)g
\]
is a Hermitian metric on~$X$ (which is K\"ahler if~$g$ is
K\"ahler) for which the restriction $d_h\restrict{X\sm U}$ of the
associated distance \fn $d_h(\cdot,\cdot)$ in~$X$ is complete.

\end{enumerate}
\end{lem}
\begin{rmk}
While it is convenient to have part (b) stated in the above form, 
if the condition holds for some choice of~$K$, then it holds
for any nonempty closed \con set on which $\psi$ is bounded.
\end{rmk}
\begin{pf*}{Proof of Lemma~\ref{completing metric on sublevel lem}}
For the proof of part~(a), observe that the $(1,1)$-form
associated to~$h$ is given by
\begin{align*}
\omega_h&\equiv\omega_g+i\partial\dbar(-\log(-\psi))
=\omega_g-\psi\inv i\partial\dbar\psi
+\psi^{-2}i\partial\psi\wedge\dbar\psi\\
&\geq (1-\delta)\omega_g+\psi^{-2}i\partial\psi\wedge\dbar\psi
\geq(1-\epsilon)\omega_g.
\end{align*}
Given two points $p,q\in Y$ and a \cinf path~$\gamma$ in~$Y$ from
$p$ to $q$, we have
\[
\ell_h(\gamma)=\int_0^1|\dot\gamma(t)|_h\,dt\geq\int_0^1\frac
1{\sqrt 2}\frac{|d\psi(\dot\gamma(t))|}{(-\psi(\gamma(t)))}\,dt.
\]
Hence \(\text{dist}_h(p,q)\geq(1-\epsilon)\text{dist}_g(p,q)\) and
\(\text{dist}_h(p,q) \geq\frac 1{\sqrt
2}\left|\log\left(\psi(q)/\psi(p)\right)\right|\), and the claim
follows. The proof of part~(b) is similar.
\end{pf*}

\begin{lem}[cf.~\cite{Nk} and \cite{Dem}]\label{completing metric on sublevel ends lem}
Suppose $(X,g,\omega_g)$ is a \con non\cpt K\"ahler manifold with
distance \fn $d_g(\cdot,\cdot)$, $E$ is an end of~$X$, $U\subset
X$ is an open set, $K$ is a nonempty closed set that contains
$X\sm E$, and $\delta\colon X\to(0,1)$ is a \cont \fnns.
\begin{enumerate}

\item[(a)] If $\vphi$ is a \cont \plsh on~$X$, $a$ is a constant
with \(\sup_E\vphi>a>\sup_{K\cap\overline E}\vphi\), and $Y$ is a
\concomp of $(X\sm E)\cup\setof{x\in E}{\vphi(x)<a}$ for which
$Y\supset K$ and $d_g\restrict{\overline Y\sm U}$ is complete,
then there exists a K\"ahler metric $h$ with distance
\fnns~$d_h(\cdot,\cdot)$ on~$Y$ \st $h\geq(1-\delta)g$ on~$Y$,
$h=g$ at each point in~$K$, and $d_h\restrict{Y\sm U}$ is
complete.

\item[(b)] If $K\cap\overline E$ is \cptns,  $d_g\restrict{X\sm(E\cup U)}$ is complete,  and there exists a \cont \plsh \fnns~$\vphi$ on~$X$ that exhausts
$\overline E$, then
there exists a K\"ahler metric $h$ with distance
\fnns~$d_h(\cdot,\cdot)$ on~$X$ \st $h\geq(1-\delta)g$ on~$X$,
$h=g$ at each point in~$K$, and $d_h\restrict{\overline E\cup(X\sm
U)}$ is complete.

\end{enumerate}
\end{lem}
\begin{pf}
For the proof of part~(a), let us fix a constant~$b$ with
$\sup_E\vphi>a>b>\sup_{K\cap\overline E}\vphi$. Then the
\fnns~$\rho$ defined by
\[
\rho\equiv
\begin{cases}
\max(\vphi-a,b-a)&\text{on }Y\cap\overline E\\
b-a&\text{on }K
\end{cases}
\]
is a negative \cont \plsh \fn on~$Y$ that approaches~$0$
at~$\partial Y$ and is equal to~$b-a$ on a \nbd of~$K$. Fixing a
\cont \fn $\eta$ with $\rho<\eta<0$ on~$Y$, the \cinf
approximation theorem of Greene and Wu (Theorem~\ref{cinf
approximation thm}) then provides a \cinf \fnns~$\psi$ on~$Y$ \st
$\psi\equiv b-a$ on a \nbd of~$K$, $\rho\leq\psi<\eta<0$ on~$Y$,
and \(\lev\psi\geq\frac 12\delta\eta g\geq\frac 12\delta\psi g\).
By Lemma~\ref{completing metric on sublevel lem}, the Hermitian
metric \(h\equiv g+\lev{-\log(-\psi)}\) has the required
properties. The proof of part~(b) is similar.
\end{pf}

We will actually prove a more general version of Theorem~\ref{cup product cont very simple version intro thm} for
a manifold with multiple ends.
For this, it will be convenient to first recall some terminology.
\begin{defn}\label{Bdd geom along set defn}
For $S\subset X$ and $k$ a nonnegative integer, we will say that a
Hermitian manifold~$(X,g)$ of dimension~$n$ {\it has bounded
geometry of order $k$ along $S$} if for some constant $C>0$ and
for every point $p\in S$, there is a biholomorphism $\Psi$ of the
unit ball $B\equiv B_{g_{\C^n}}(0;1)\subset\C^n$ onto a \nbd of
$p$ in $X$ \st $\Psi(0)=p$ and \st on $B$,
\[
C\inv g_{\C^n}\leq\Psi^*g\leq Cg_{\C^n}
\quad\text{and}\quad|D^m\Psi^*g|\leq C\text{ for }m=0,1,2,\dots,k.
\]
\end{defn}

\begin{defn}[cf.~Definition~2.2
of~\cite{NR Filtered ends}]\label{Special end definition} We will
call an end $E$ of a \con non\cpt complete Hermitian manifold $X$
\emph{special} if $E$ is of at least one of the following types:
\begin{enumerate}
\item[(BG)] $X$ has bounded geometry of order $2$ along $E$;

\item[(W)] There exists a \cont \plsh \fn $\vphi$ on~$X$ \st
\[
\setof{x\in E}{\vphi(x)<a}\Subset X\qquad\forall\,a\in\R;
\]

\item[(RH)] $E$ is a regular hyperbolic end (i.e., $E$ is a hyperbolic end 
and the Green's \fn on~$X$ vanishes
at infinity along $E$); or

\item[(SP)] $E$ is a parabolic end, the Ricci curvature of $g$ is
bounded below on $E$, and there exist positive constants $R$ and
$\delta$ \st \(\Vol\,\big(B(x;R)\big)>\delta\) for all $x\in E$.

\end{enumerate}
We will call an ends decomposition in which each of the ends is special a \emph{special ends decomposition}.
\end{defn}

While bounded geometry and special ends do not play fundamental roles in this paper, these conditions have been shown to strongly
determine the holomorphic structure of complete K\"ahler manifolds.
In particular, according to \cite{Gro-Sur la groupe fond}, \cite{Li Structure
complete Kahler}, \cite{Gro-Kahler hyperbolicity},
\cite{Gromov-Schoen}, \cite{NR-Structure theorems},
\cite{Delzant-Gromov Cuts}, \cite{NR Filtered ends}, and \cite{NR L2 Castelnuovo}, a \con non\cpt complete K\"ahler manifold that admits
a special ends decomposition and has at least three (filtered) ends admits a proper \holo mapping onto a
Riemann surface. 
Applications of the versions
of the cup product lemma from 
Section~\ref{cup product cont sect} in the above context will be considered elsewhere (see~\cite{NR Weakly Special}).

Theorem~\ref{BH plsh fn cont two pos ev thm} is an immediate
consequence of the following
 (cf.~Proposition~4.4 of~\cite{NR Weak Lefschetz}):
\begin{thm}\label{plsh fn two pos ev par ends thm}
Let $(X,g)$ be a \con non\cpt complete K\"ahler manifold that
admits a \cont \plsh \fnns~$\vphi$ whose restriction to some 
$2$-dimensional germ of an \anal set 
at some point is \str \plshns.
\begin{enumerate}

\item[(a)] For every compactly supported $\dbar$-closed
\cinf form~$\alpha$ of type~$(0,1)$ on~$X$, there exists a bounded
\cinf \fnns~$\beta$ on~$X$ \st $\dbar\beta=\alpha$ on~$X$ and
$\beta$ vanishes on every end $E$ of~$X$ for which
$E\cap\supp\alpha=\emptyset$ and $\sup_E\vphi>\sup_{\partial
E}\vphi$ (here, we take $\sup\emptyset=-\infty$).

\item[(b)] In any ends decomposition $X\sm K=E_1\cup\cdots\cup
E_m$ with $K\neq\emptyset$, $\sup_{E_j}\vphi>\sup_{\partial
E_j}\vphi$ for exactly one choice of~$j$, say $j=1$,  the
remaining ends $E_2,\dots,E_m$ are parabolic and not special, and
every \holo \fn on~$E_1$ admits a (unique) extension to a \holo
\fn on~$X$. Moreover, for every nonconstant \cont \plsh \fnns~$\psi$ on~$X$, 
$\sup_{E_1}\psi>\sup_{\partial E_1}\psi$
while $\sup_{E_j}\psi=\sup_{\partial E_j}\psi$ for $j=2,\dots,m$.

\item[(c)] If $e(X)=1$ (i.e., $\sup_E\vphi>\sup_{\partial E}\vphi$
for every end~$E$), then $H^1_c(X,\ol)=0$.

\end{enumerate}
\end{thm}
\begin{pf}
Let $n\equiv\dim X$. It follows from the condition on~$\vphi$
that there exists a $2$-dimensional \con \cpx submanifold $Z$
of a nonempty \rel \cpt domain~$U$ in~$X$ \st
$\vphi\restrict Z$ is \str \plshns, and by replacing $\vphi$ with $e^\vphi$
if necessary, we may assume that $\vphi>0$ on~$X$. Given an ends
decomposition $X\sm K=E_1\cup\cdots\cup E_m$ for~$X$ with
$K\neq\emptyset$, which we may order so that
$\sup_{E_1}\vphi=\sup\vphi$, part~(a) of Lemma~\ref{basic ends
lem} provides a \con \cpt set~$K'$ \st $K\cup U\subset K'$ and any
domain $Y\supset K'$ has an ends decomposition $Y\sm
K=E_1'\cup\cdots\cup E_m'$, where $E_j'=E_j\cap Y$ for
$j=1,\dots,m$. Fixing a constant~$a$ with
$\max_{K'}\vphi<a<\sup\vphi$, we may let the above domain~$Y$ be
the \concomp of~$\setof{x\in X}{\vphi(x)<a}$ containing~$K'$.

Lemma~\ref{completing metric on sublevel ends lem} and
Theorem~\ref{cup product cont very simple version intro thm} (or
Theorem~\ref{cup product cont simple version thm}) together provide
a complete K\"ahler metric~$h$ on~$Y$ \st
$\partial\vphi\wedge\theta\equiv 0$ for every $L^2$ (closed) \holo
$1$-form~$\theta$ on $(Y,h)$. However, if $f$ is a nonconstant
\holo \fn on a nonempty domain $W\subset U$, $p\in W\cap Z$, and
$L$ is the level of~$f$ through~$p$, then $L\cap Z$ is an \anal
set of positive dimension at~$p$ and $\vphi\restrict{(L\cap Z)}$ is
\str \plshns. It follows from Lemma~\ref{plsh fn holo 1 form wedge 0 lem} that
$\partial\vphi\wedge df$ is not everywhere~$0$ in~$W$, and hence $(Y,h)$ has no nontrivial $L^2$ \holo $1$-forms.

Applying part~(b) of Theorem~\ref{Bochner Hartogs no L2 forms
thm}, we see that for $j=2,\dots,m$, $E_j'$ is a parabolic end of
$(Y,h)$, and hence $E_j=E_j'$ and
$\sup_{E_j}\vphi\leq\max_{\partial E_j}\vphi$. 
If $f\in\ol(E_1)$, then there exists a \holo \fnns~$u_0$
on $Y\supset (X\sm E_1)\cup K'$ \st $u_0=f$ on $E_1'=E_1\cap Y$,
and hence the \fn $u$ on~$X$ given by $u\equiv f$ on~$E_1$ and
$u\equiv u_0$ on~$Y$ is a \holo extension of~$f$ to~$X$.  Furthermore, if 
there exists a \cont \plsh \fnns~$\psi$ on~$X$ with $\sup_{E_j}\psi>\max_{\partial E_j}\psi$
for some $j>1$, then after replacing~$\psi$ with the composition of a suitable nondecreasing \cvx \fn and~$\psi$,
we may assume that $\sup_{E_j}\psi=\infty$. But then for $\epsilon>0$ sufficiently small, the \plsh
\fn $\vphi_0\equiv\vphi+\epsilon\psi$, which satisfies the condition placed on~$\vphi$ in the statement
of the theorem, must also satisfy $\sup_{E_1}\vphi_0>\max_{\partial E_1}\vphi_0$ and
 $\sup_{E_j}\vphi_0=\infty=\sup\vphi_0$, which as we have seen, is impossible. Thus such a 
\fnns~$\psi$ cannot exist.

Part (b) will follow if we prove that for each $j=2,\dots,m$, with
respect to~$g$ in~$X$, $E_j$ is neither a hyperbolic end nor a
special end. For this, let us first fix an ends decomposition
$X\sm C=A_0\cup A_1$ \st $A_0\subset E_j$, $E_j\sm A_0\Subset X$,
and $A_1\supset(E_1\cup\cdots\cup\what{E_j}\cup\cdots\cup E_m)$,
as provided by part~(d) of Lemma~\ref{basic ends lem}.  We may
also choose the above domain $Y$ so large that $C\subset Y$ and we
have the ends decomposition $Y\sm C=A_0\cup A_1'$, where
$A_1'=A_1\cap Y$ (and $A_0=A_0\cap Y$). Applying Lemma~\ref{completing
metric on sublevel ends lem}, we get a complete K\"ahler metric
$k$ on~$Y$ \st $k=g$ at each point in the closed set $(X\sm
E_1)\cup K'\supset A_0$. In particular, $A_1'$ is a hyperbolic end
in $(Y,k)$, and if $E_j$ is a hyperbolic end or a special end of
type~(SP) with respect to~$g$ (the latter holding
if, for example, $E_j$ is parabolic and of type~(BG)), then $A_0$ must be hyperbolic or of
type (SP) with respect to $k$, and by applying Theorem~3.6
of~\cite{NR L2 Castelnuovo} in~$(Y,k)$, we are able to produce a
nonnegative nonconstant \cinf \plsh \fn on~$Y$ that vanishes on
$Y\sm E_j$ and in particular, extends by~$0$ to a \cinf \plsh
\fnns~$\psi$ on~$X$ with $\sup_{E_j}\psi>0=\max_{\partial E_j}\psi$. 
But as we have seen, such a \fnns~$\psi$ cannot exist.  Similarly, $E_j$ cannot be of type (W).
Thus part~(b)
is proved.

For the proof of part~(a), letting $\alpha$ be a nontrivial
compactly supported $\dbar$-closed \cinf form of type
$(0,1)$ on~$X$, we may choose the above \cpt set~$K$ to
contain~$\supp\alpha$. Part~(a) of Theorem~\ref{Bochner Hartogs no
L2 forms thm} then provides a bounded 
\cinf \fnns~$\beta_1$ with finite
$h$-energy \st $\dbar\beta_1=\alpha$ on~$Y\supset K$ and
$\beta_1\equiv 0$ on any hyperbolic end of $Y$ that does not meet
$\supp\alpha$. In particular, $\beta_1\equiv 0$ on $E_1'$, so
$\beta_1$ extends to a bounded \cinf \fnns~$\beta$ on $X=Y\cup E_1$ that
vanishes on~$E_1$. Furthermore, if $E\subset X\sm\supp\alpha$ is
an end of $X$ and $\sup_E\vphi>\sup_{\partial E}\vphi$, then there
is end $E_0$ of $X$ \st $E_0\subset E_j\cap E$ for some $j$ and
$\sup_{E_0}\vphi>\sup_{\partial E_0}\vphi$. On the one hand, since
$E_l$ is a parabolic end and $\vphi$ is bounded on~$E_l$ for each
$l=2,\dots,m$, we must have $j=1$; i.e., $E_0\subset E_1\cap E$.
On the other hand, by the above, there exists a bounded \cinf
\fnns~$\beta'$ on~$X$ \st $\dbar\beta'=\alpha$ and $\beta'\equiv
0$ on $E$. The \holo \fn $\beta-\beta'$ vanishes on $E_0$, and
therefore on~$X$, so we must have $\beta=\beta'$, and hence
$\beta\equiv 0$ on~$E$.

Part~(c) follows from parts~(a) and (b).
\end{pf}

\bibliographystyle{amsalpha.bst}

\end{document}